%% file: local_timestepping_paper.tex
\documentclass[onefignum,onetabnum]{siamart171218}

\pdfoutput=1

\usepackage{amssymb,amsmath}
\usepackage{hyperref}
\usepackage{graphicx}
\usepackage{tikz}
\usepackage{cancel}
\usepackage{mathtools}
\usepackage{array}
\usepackage{hhline}
\usepackage{longtable}

\headers{High-order, conservative local time-stepping}{
  W. Throwe and S. A. Teukolsky}

\title{A high-order, conservative integrator with local time-stepping%
  \thanks{Submitted to the editors 10 October 2019.
    \funding{
      This work was supported in part by the Sherman Fairchild Foundation and
      NSF grants \hbox{PHY-1606654} and \hbox{ACI-1713678} at Cornell.}}}

\newcommand\Cornell{Center for Astrophysics and Planetary Science, Cornell
  University, Ithaca, New York 14853, USA}
\newcommand\Caltech{Theoretical Astrophysics, Walter Burke Institute for
  Theoretical Physics, California Institute of Technology, Pasadena, CA 91125,
  USA}

\author{William Throwe\thanks{\Cornell{} (\email{wtt6@cornell.edu})} \and
  Saul A. Teukolsky\thanks{\Cornell{} and \Caltech}}

\usepackage{soul}
\setstcolor{red}

\makeatletter
\newcommand\Biggg[1]{{\hbox{$\left#1\vbox to 23.5pt{}\right.\n@space$}}}
\makeatother

\newcommand\abs[1]{\lvert#1\rvert}

\newcommand\deriv[2][]{\frac{d#1}{d#2}}
\newcommand\pderiv[2][]{\frac{\partial#1}{\partial#2}}
\newcommand\derivn[3][]{\frac{d^{#3}#1}{d#2^{#3}}}

\newcommand\bof[1]{\!\big(#1\big)}

\newcommand\mat{\mathbf}

\newcommand\dt{\Delta t}
\newcommand\dx{\Delta x}
\newcommand\y{\mat{y}}
\newcommand\dy{\Delta \y}
\newcommand\D{\mat{D}}

\newcommand\ty{\tilde{\y}}
\newcommand\dty{\Delta \tilde{\y}}
\newcommand\tU{\tilde{t}}
\newcommand\dtU{\Delta \tilde{t}}
\newcommand\tD{\tilde{\D}}

\newcommand\V{\mat{V}}
\newcommand\B{\mat{B}}

\newcommand\allq{q^1\cdots q^S}
\newcommand\ofyq{\big(\y^1_{q^1}, \ldots, \y^S_{q^S}\big)}

\newcommand\charspeed{\lambda}

\newcommand\etal{~{\it et~al.}}

\newcommand\Astep[3][]{\draw[#1] (-1/8, #2) -- (-3/4, #2) node[left] {%
    \setbox0\hbox{$-2\dt^A$}%
    \hbox to\wd0{\hfil$#3$\hfil}};}%
\newcommand\Bstep[3][]{\draw[#1] (1/8, #2) -- (3/4, #2) node[right] {%
    \setbox0\hbox{$-2\dt^A$}%
    \hbox to\wd0{\hfil$#3$\hfil}};}%
\newcommand\Astepunused[1]{\Astep[dotted]{#1}{}}
\newcommand\Bstepunused[1]{\Bstep[dotted]{#1}{}}
\newcommand\Achart[2]{\node at (-7/16, #1) {#2};}
\newcommand\Bchart[2]{\node at (7/16, #1) {#2};}
\newcommand\chartdots[1]{\node at (0, #1) {$\vdots$};}

\newcolumntype{M}{>{$}c<{$}}
\newcolumntype{s}{>{\small$}c<{$}}

\begin{document}

\maketitle

\begin{abstract}
  We present a family of multistep integrators based on the Adams-Bashforth
  methods.  These schemes can be constructed for arbitrary convergence order
  with arbitrary step size variation.  The step size can differ between
  different subdomains of the system.  It can also change with time within a
  given subdomain.  The methods are linearly conservative, preserving a wide
  class of analytically constant quantities to numerical roundoff, even when
  numerical truncation error is significantly higher.  These methods are
  intended for use in solving conservative PDEs in discontinuous Galerkin
  formulations
  or in finite-difference methods with compact stencils.
  A numerical test
  demonstrates these properties and shows that significant speed improvements
  over the standard Adams-Bashforth schemes can be obtained.
\end{abstract}



\section{Introduction}

A common problem in computational fields is to find approximate solutions to
partial differential equations~(PDEs).  For hyperbolic PDEs, where a solution
typically describes an evolution of one or more fields through time, the most
common approach is to apply the method of lines, where the spatial coordinates
in the PDE are discretized, producing a large system of coupled ordinary
differential equations~(ODEs)
with one degree of freedom per variable per grid point.
These systems of equations can then be
discretized in time and solved using standard explicit integration schemes.

In order to obtain an accurate solution, it is a necessary condition that
the time discretization must be
fine enough for the integration to be stable.  For a method-of-lines system,
this limit is primarily because of the Courant-Friedrichs-Lewy~(CFL) condition,
which limits the step size to approximately the information propagation time
between grid points.  The resulting step size can show large variation across
the spatial domain because of changes in the propagation speed or, more
commonly, because of changes in the spacing of the evaluation points.  It is
often desirable to increase the density of points in some regions to
resolve
rapidly varying parts of the solution, but this then restricts the step size
allowed for stability.  Furthermore, in order to evaluate the system right-hand
side, it is necessary to know the entire state of the system at the time of
interest.  The time step for the whole system is then set by the most
restrictive of the conditions over the entire domain.  If the problematic
points make up a small fraction of the system, then the forced evaluations at
the remaining points can dominate the computational expense.

To reduce the computational cost of finding these solutions, we would like to
evaluate each point at intervals set by its own stability limit, rather than
the smallest limit for all the points.  A method allowing this is known as a
\emph{local time-stepping}~(LTS) (or \emph{multirate}) method, as opposed to a
global time-stepping~(GTS) method.  Such a method must describe an update
scheme for the frequently evaluated degrees of freedom that does not require
knowing the full state of the system.

Modifying a GTS method into an LTS one can have significant drawbacks.  The
individual steps near locations of time step changes are typically more
expensive than for a GTS method, so the benefit of fewer derivative evaluations
must outweigh this overhead.  Care must be taken when calculating the CFL limit
near step size changes to take into account variations in the characteristic
speeds of the system in the neighborhood of the element.~\cite{Gnedin2018}
Furthermore, modifying the GTS scheme can destroy numerically desirable
properties of the integrator, such as a high convergence order.  LTS schemes
also do not naturally provide exact conservation of linear conserved
quantities~\cite{Sandu2009}, a property
usually guaranteed by
GTS integrators.
In a physical system, errors accumulated in these quantities
(which can represent, for example, total mass) can produce an approximate
solution qualitatively different from the true solution.

Early LTS schemes (for example~\cite{Berger1984,Gear1984}) typically used GTS
integrators with different time steps and performed interpolation to obtain
data at times at which it was not produced directly.  Such schemes are easy to
adapt to arbitrary mesh configurations and can be constructed to obtain the
same convergence order as the underlying GTS method, but they do not preserve
conserved quantities of the system.  Corrections to more accurately treat
conservation laws were developed~\cite{Berger1987}, but still only resulted in
approximate conservation.

More recently, many methods have been investigated as starting points for more
sophisticated LTS methods, including both substep~\cite{Krivodonova2010,
  Gassner2011, Grote2015, Fok2016, Gunther2016, Almquist2017} and
multistep~\cite{Sandu2009, Grote2013, Winters2014} integrators and also less
common methods such as leapfrog~\cite{Grote2013, Grote2018}, Richardson
extrapolation~\cite{Constantinescu2013}, ADER~\cite{Schoeder2018}, and implicit
methods~\cite{Chabaud2012}.  Demirel\etal~\cite{Demirel2015} have even explored
LTS schemes constructed from multiple unrelated GTS integrators.  Recently,
G\"unther and Sandu~\cite{Gunther2016} presented a very general family of
multirate Runge-Kutta-like methods based on the GARK family of
integrators~\cite{Sandu2015} that unifies many of the previous
Runge-Kutta-based LTS schemes.  These methods are applicable to any problem and
can be constructed to have any order of accuracy, but they are not
conservative.  Sandu and Constantinescu~\cite{Sandu2009} presented an
Adams-Bashforth-based scheme based on evaluating the right-hand side of the
evolution equations using a combination of data at different times.  This
system is conservative and applicable to any system of equations, but the
method is limited to second-order accuracy at times at which all degrees of
freedom are evaluated and first-order accuracy at intermediate times.

LTS integrators for the special case of linear systems have been developed
based on Adams-Bashforth~\cite{Grote2013}, Runge-Kutta~\cite{Grote2015,
  Almquist2017}, and leapfrog~\cite{Grote2013, Grote2018} schemes.  Of
particular interest here, starting from the Adams-Bashforth methods, Grote and
Mitkova~\cite{Grote2013} found a family of high-order, conservative methods for
integer ratios between step sizes on different degrees of freedom.  These
methods use the linearity of the system to split the equations into a form
resembling multiple copies of the standard Adams-Bashforth method.

Some authors have derived methods specialized to the discontinuous Galerkin or
finite volume formalisms.  The structure of elements coupled comparatively
weakly in a standard way by exchange of fluxes allows for some simplifications
to the problem.  Winters and Kopriva~\cite{Winters2014} presented a scheme
using dense output of the integrators for each element to calculate fluxes at
intermediate times.  This scheme is high-order and allows for arbitrary step
ratios and varying time steps, but it sacrifices the conservative nature of its
parent scheme.  Gassner\etal~\cite{Gassner2011} presented a similar method, but
restored conservation by treating the element and flux terms as a predictor and
corrector.  Krivodonova~\cite{Krivodonova2010} constructed a method based on a
Runge-Kutta integrator which, while not naturally conservative, was made so by
adding a correction to cancel any error in conservation whenever neighboring
cells are aligned in time.  Cavalcanti\etal~\cite{Cavalcanti2015} considered
the addition of nonlinear operations, such as slope limiting, to the
integration step.

In this paper we present a high-order, conservative
scheme based on the Adams-Bashforth family of explicit multistep methods
capable of solving any explicit initial value problem.
The method uses the
idea of performing single right-hand side evaluations using values from
different times, in a similar manner to previous work presented by Sandu and
Constantinescu~\cite{Sandu2009}.  The scheme is conservative and has the same
convergence order as the Adams-Bashforth integrator it is based on.  The method
allows for generic ratios of step sizes between different degrees of freedom,
as well as for arbitrarily varying the time steps of the individual degrees of
freedom.
The applications discussed here are to discontinuous Galerkin methods, but the
method can also be applied to finite-difference schemes with compact stencils.
(Although the method is in principle fully general for any explicit initial
value problem for a set of coupled ODEs, in practice it becomes extremely
inefficient if the spatial couplings are non-local.)
When applied to a linear system with integer step size ratios,
this scheme reduces to the Adams-Bashforth-based scheme presented by Grote and
Mitkova~\cite{Grote2013}.

Throughout this work, we refer to approximations with leading order error
proportional to $\dt^{n+1}$ as ``order-$n$'' approximations or as ``accurate
to order
$n$.''  Thus, a single step of order-$k$ Adams-Bashforth is an order-$k$
approximation of the new value.

The remainder of this paper is structured as follows: Section~\ref{sec:method}
presents a derivation of the integration scheme.
Section~\ref{sec:special-cases} discusses simplifications that are applicable
when the method is applied to some common special cases.
Section~\ref{sec:numerical} applies the method to numerical test cases.
Appendices~\ref{sec:consistency} and \ref{sec:error} demonstrate the
consistency of the method and compute the leading order error term.
Appendix~\ref{sec:general-splitting} provides a method of writing common GTS
integrators in a form mathematically similar to the LTS integrator presented
here.  Finally, Appendix~\ref{sec:tables}
lists specific formulas for methods of order 2, 3, and 4.

\section{The method}
\label{sec:method}

\subsection{Adams-Bashforth methods}

Suppose we wish to numerically solve a set of coupled first-order ordinary
differential equations
\begin{equation}\label{eq:ode}
  \deriv[\y]{t} = \D(\y),
\end{equation}
with $\y = \y_0$ at some initial time $t_0$.  Here
$\D(\y)$, the time-derivative operator, is the right-hand side evaluated
when the system is in state $\y$.
(Throughout this work we will restrict our attention to autonomous systems.  As
any non-autonomous system can be recast in an autonomous form (see, for
example, Section II.2 of \cite{Hairer:ODE1}) we lose no generality from this
restriction.)
A common method is to solve for the
variables at a (monotonic) sequence of times $t_0, t_1, \ldots$ using a
$k$th-order Adams-Bashforth method
\begin{equation}\label{eq:gts}
  \dy_n = \dt_n \sum_{j=0}^{k-1} \alpha_{nj} \D(\y_{n-j})
\end{equation}
with $\dy_n = \y_{n+1} - \y_n$, $\dt_n = t_{n+1} - t_n$ and the coefficients
corresponding to the step given by~\cite{Winters2014}
\begin{equation}
  \alpha_{nj} = \frac{1}{\dt_n} \int_{t_n}^{t_{n+1}} dt\,
  \ell_j\!\!\left(t; t_n, t_{n-1}, \ldots, t_{n-(k-1)}\right).
\end{equation}
Here
\begin{equation}
  \ell_n(t; t_0, \ldots, t_{k-1}) =
  \prod_{\substack{j=0\\j\ne n}}^{k-1} \frac{t - t_j}{t_n - t_j}
\end{equation}
are Lagrange polynomials.
The values required for evaluating the first step can be obtained using a
standard Adams-Bashforth GTS self-start procedure~\cite{Mersman1965}, after
which the step sizes can be adjusted to the desired ratio.

If different degrees of freedom require different time steps for stability, it
may be desirable to evaluate these variables at different frequencies, in order
to avoid unnecessary computations for the more stable variables.
Suppose we divide the components of $\y$ into $S$ sets $\y^1, \ldots, \y^S$,
each of which is a collection of degrees of freedom that are always evaluated
at the same times,
$t^s_0, t^s_1, \ldots$.
We can then split~\eqref{eq:ode} into
an equation for each of these sets:
\begin{equation}
  \deriv[\y^s]{t} = \D^s\!\!\left(\y^1, \ldots, \y^S\right),
\end{equation}
where $\D^s$ is the result of $\D$ restricted to the set~$s$.  Any attempt to
use this equation to perform an LTS evolution immediately encounters the
problem that evaluating its right-hand side requires knowing the entire state
of the system, which conflicts with the goal of independent evaluation times
for different degrees of freedom.

\subsection{Conserved quantities}

A linear conserved quantity is a quantity~$C$ expressible as an inner product
of a vector~$\mat{c}$ with the evolved variables (treated as a vector)
\begin{equation}\label{eq:conserved}
  C = \mat{c} \cdot \y,
\end{equation}
with
\begin{equation}\label{eq:conserved-derivative}
  \mat{c} \cdot \D(\y) = 0
\end{equation}
for all values of $\y$.  Such a quantity is constant under exact integration of
the system and under integration using Euler's method.  An integrator is called
(linearly) \emph{conservative} if all such quantities remain constant when
integrating a system using it~\cite{Shampine1986-conservation}.  It is
desirable for an integrator to keep such quantities precisely constant (up to
roundoff error) rather than merely constant up to the truncation error of the
scheme.  Such quantities often have an intuitive physical meaning, and
frequently even a small rate of drift can cause qualitative changes in the
evolution of the system.

When solving a PDE representing a physical system, the most common linear
conserved quantities are integrals over the computational domain of fields
representing densities.  The vector $\mat{c}$ in these cases is the
vector of coefficients necessary to perform a numerical integral.  In a
discontinuous Galerkin scheme these coefficients would combine quadrature
weights on the elements and factors arising from coordinate mappings of the
elements relative to their canonical shapes.

\subsection{Second-order $2:1$ stepping}

\begin{figure}
  \centering%
  \begin{tikzpicture}
    \Astepunused{-2}
    \Astep{-1}{-\dt^A}
    \Astep{0}{0}
    \Astep{1}{\dt^A}
    \Bstepunused{-4/2}
    \Bstepunused{-3/2}
    \Bstepunused{-2/2}
    \Bstep{-1/2}{-\dt^B}
    \Bstep{0}{0}
    \Bstep{1/2}{\dt^B}
    \Bstep{2/2}{2\dt^B}
    \Achart{1/2}{(a)}
    \Bchart{1/4}{(b)}
    \Bchart{3/4}{(c)}
  \end{tikzpicture}
  \caption{The step pattern for a $2:1$ method on two sets, with time steps
    $\dt^A = 2 \dt^B$.  There are three types of steps: the large step on
    element A marked (a), and the two types of small step on B marked (b) and
    (c).  For a second-order method, we use only the two most recent values of
    the variables when taking a step.  Steps whose values are no longer needed
    for the indicated steps are marked with dotted lines.}
  \label{fig:initial-example}
\end{figure}
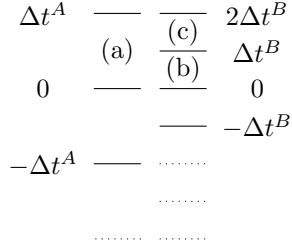

Let us first consider as an example the case of a second-order scheme on two
sets, $A$ and $B$, with $B$ being evaluated twice as often as $A$.  Call their
step sizes~$\dt^A$ and $\dt^B = \dt^A/2$.  This step pattern is shown in
Figure~\ref{fig:initial-example}, where for simplicity we consider the steps
starting from~$t=0$ leading up to~$t=\dt^A$.  There are three types of steps to
consider: the large step on set~$A$, labeled $(a)$, and the first and second
halves of that step on set~$B$, labeled $(b)$ and $(c)$.  This case is
considered in Sandu and Constantinescu~\cite{Sandu2009}, but the method
presented there only provides a second-order value when sets have stepped to
the same time; intermediate values are only accurate to first order.

We will start with the small step~$(b)$.  For a GTS Adams-Bashforth method,
this step would be given by
\begin{equation}\label{eq:example-b-tilde}
  \dy^B_b = \dt^B \left[
    \frac{3}{2} \tD^B\bof{0} - \frac{1}{2} \tD^B\bof{-\dt^B}
    \right],
\end{equation}
where
$\tD^B(t)$ is an approximation to the derivative at time~$t$ that we will
now construct.

At time~$t=0$ we have values for the entire system, so
we can obtain a derivative by applying the derivative operator to the initial
data, so we take
$\tD^B(0) = \D^B(\y^A(0), \y^B(0))$.
We cannot evaluate $\tD^B(-\dt^B)$ in this manner, however,
because we do not have data
for $\y^A$ at $t=-\dt^B$, so we must construct it from the values at $t=0$ and
$t=-2\dt^B$.
Up through this point, the derivation matches the results of Sandu and
Constantinescu~\cite{Sandu2009}, but they now choose to use the known value of
$\y^A(-2\dt^B)$ to evaluate the derivative, while we will search for a more
accurate approximation.
There are two reasonable choices of how to do this: average the
known values of $\y^A$ to get a value at the desired time and use that to apply
the derivative operator, or apply the derivative operator at both times (using
the value of $\y^B$ at $-\dt^B$ both times) and average the results.  We choose
the latter,
taking $\tD^B(-\dt^B) = [\D^B(\y^A(0), \y^B(-\dt^B)) + \D^B(\y^A(-\dt^A),
\y^B(-\dt^B))] / 2$, so \eqref{eq:example-b-tilde} becomes
\begin{multline}\label{eq:example-b}
  \dy^B_b = \dt^B \bigg[
    \frac{3}{2} \D^B\bof{\y^A\bof{0}, \y^B\bof{0}}
    - \frac{1}{4} \D^B\bof{\y^A\bof{0}, \y^B\bof{-\dt^B}} \\
    - \frac{1}{4} \D^B\bof{\y^A\bof{-\dt^A}, \y^B\bof{-\dt^B}}
    \bigg].
\end{multline}
The error in averaging the derivatives is of order~$(\dt^B)^2$, so it
introduces an error of order~$(\dt^B)^3$ in the value after the step,
preserving the second-order quality of the base GTS method.

The second small step, $(c)$, proceeds similarly, except that we now use a
derivative at $\dt^B$ instead of $-\dt^B$.  Instead of averaging the
derivatives at different $\y^A$ we must therefore perform a (linear)
extrapolation to obtain our approximate
derivative~$\tD^B(\dt^B) = [3\, \D^B(\y^A(0), \y^B(\dt^B)) - \D^B(\y^A(-\dt^A),
\y^B(\dt^B))] / 2$.  Thus,
we obtain the rule
\begin{multline}\label{eq:example-c}
  \dy^B_c = \dt^B \bigg[
    \frac{9}{4} \D^B\bof{\y^A\bof{0}, \y^B\bof{\dt^B}}
    - \frac{3}{4} \D^B\bof{\y^A\bof{-\dt^A}, \y^B\bof{\dt^B}} \\
    - \frac{1}{2} \D^B\bof{\y^A\bof{0}, \y^B\bof{0}}
    \bigg].
\end{multline}

We could use the same procedure to evaluate the large step~$(a)$, but, as this
would not take into account the value $\y^B(\dt^B)$ used for taking the second
small step, there is no way this procedure could be conservative.  This,
however, gives us a hint as to how to proceed: we treat the large step as
having two internal steps, one for balancing each of the small steps.  In fact,
in order to remain conservative, we must take each of these internal steps
using the same scheme as for the corresponding small step, except using the
part of the derivative corresponding to set~$A$.  This can be seen by
considering a generic pair of methods for a step:
$\dy^{A,B} = \sum_i k^{A,B}_i \D^{A,B}(\mat{q}_i)$, where the $k^{A,B}_i$ are
given coefficients and the $\mat{q}_i$ are vectors obtained in some manner from
the known values of $\y$.
The change in a linear conserved quantity during that step is
\begin{align}
  \Delta C^A + \Delta C^B
  &= \mat{c}^A \cdot \sum_i k^A_i \D^A(\mat{q}_i)
   + \mat{c}^B \cdot \sum_i k^B_i \D^B(\mat{q}_i) \\
  &= \sum_i \left[k^A_i \mat{c} \cdot \D(\mat{q}_i)
    + (k^B_i - k^A_i) \mat{c}^B \cdot \D^B(\mat{q}_i)\right].
\end{align}
The first term vanishes by~\eqref{eq:conserved-derivative}, so the only way for
two sets to take equal-sized steps in a conservative manner is if they use the
same step rule.  The procedure for the large step can therefore be found by
summing~\eqref{eq:example-b} and \eqref{eq:example-c}, giving
\begin{multline}
  \dy^A_a = \\ \dt^A \bigg[
    \frac{9}{8} \D^A\bof{\y^A\bof{0}, \y^B\bof{\dt^B}}
    + \frac{1}{2} \D^A\bof{\y^A\bof{0}, \y^B\bof{0}}
    - \frac{1}{8} \D^A\bof{\y^A\bof{0}, \y^B\bof{-\dt^B}} \\
    - \frac{3}{8} \D^A\bof{\y^A\bof{-\dt^A}, \y^B\bof{\dt^B}}
    - \frac{1}{8} \D^A\bof{\y^A\bof{-\dt^A}, \y^B\bof{-\dt^B}}
    \bigg].
\end{multline}
Note that the coefficients have changed by a factor of $2$ compared to the
previous equations because of the change of the leading coefficient to $\dt^A$.
As the two small steps were accurate to second order and this is effectively
their concatenation, it is also accurate to second order.
A detailed analysis of the consistency and convergence of the general method
described below, of which this is a special case,
is presented in Appendix~\ref{sec:consistency}.

\subsection{Conservative time steppers}

Let us turn now to the task of finding a general conservative, high-order LTS
integrator.  First, we will consider the implications of requiring an
Adams-Bashforth-like LTS scheme to be conservative.  For such a scheme it only
makes sense to evaluate~\eqref{eq:conserved} at times at which all the degrees
of freedom are evaluated.  We therefore introduce a new quantity~$\ty_n$ that
is defined for the entire set of degrees of freedom for each time~$\tU_n$ at
which any set is evaluated, and is equal to
the (as yet unknown) complete numerical solution at all grid points wherever it
exists.
If we provide
a rule to compute $\ty_{n+1}$ from earlier values
then, as long as portions of $\ty_n$
that do not correspond to
values in one of the $\y^s_m$
are never used, we can obtain an LTS method by
summing the $\dty_n$
between evaluations.
Furthermore, if the step from $\ty_n$ to $\ty_{n+1}$ is conservative, then the
implied full method will be as well.

The condition for this small step to be conservative is
\begin{equation}\label{eq:conserved-delta}
  0 = \mat{c} \cdot \dty_n.
\end{equation}
This is satisfied if we evaluate $\dty_n$ using a standard Adams-Bashforth
method, but that would require values of $\ty_n$ that
are not included in any of the $\y^s_m$.
Comparing~\eqref{eq:conserved-derivative} and \eqref{eq:conserved-delta}, we
see that we will obtain a conservative method if we take
\begin{equation}
  \dty_n = \sum_{q^1\cdots q^S} \beta_{n;\allq} \D\ofyq,
\end{equation}
for some coefficients $\beta_{n;i^1\cdots i^S}$.
The choices of these coefficients
are not unique, but there is a natural choice.  We evaluate each step using a
standard order-$k$ Adams-Bashforth scheme, but instead of using the derivatives
of the function that we cannot evaluate, we use approximate derivatives
$\tD_n$.  As long as these are accurate to order~$k-1$, we will lose no formal
accuracy for the step.  We evaluate $\tD_n$ by treating $\D(\y^1(t_1), \ldots,
\y^S(t_S))$ as a function of the times $t_1, \ldots, t_S$ independently, and
then performing a multidimensional interpolation
from the known values in the space of the evaluation times on each set.
To obtain
the required accuracy,
we must use evaluations from
at least $k$ times from each set,
and it is natural to choose the most recent values.  The known values of $\D$
then form a lattice in the multidimensional space.  Multidimensional
interpolation from such a lattice is not unique, but a natural choice is to
perform it as a series of one-dimensional interpolations.\footnote{This freedom
  arises from the fact that the system of equations defining this interpolation
  is underdetermined for $S$ and $k$ greater than~1: we must find $k^S$ fitting
  coefficients but there are only $\binom{k+S-1}{S}$ monomial terms of degree
  less than $k$ (which are the ones relevant for an order $k$ fit).  A general
  choice of interpolation coefficients will result in an interpolating
  polynomial containing all terms of degree less than $k$ in each of the $t^s$
  individually.  We therefore have the freedom to modify the interpolation
  coefficients as long as the modification alters only terms in the
  interpolating polynomial of total degree at least $k$.  This freedom could be
  used, for example, to set certain coefficients to zero to reduce the number
  of computations required or to decrease the effect of terms where the times
  on different sets have large mismatches.

  In the case where the step size on each set is constant, the alternative sets
  of interpolation coefficients can be obtained by adding high-order products
  of discrete Chebyshev polynomials~\cite{Chebyshev:discrete-polynomials} to
  the coefficients in
  \eqref{eq:lts-coef-small}.  In
  the general case we know of no simple method to calculate alternative
  coefficients.  We have not investigated the use of such alternative
  coefficients in either of these cases.} Combining all these ideas, we have
\begin{equation}
  \dty_n =
  \dtU_n \sum_{i=0}^{k-1} \tilde{\alpha}_{ni}
  \sum_{q^1=0}^{k-1} \cdots \sum_{q^S=0}^{k-1}
  I_{ni;\allq}
  \D\!\left(\y^1_{m^1(n)-q^1}, \ldots, \y^S_{m^S(n)-q^S}\right),
\end{equation}
where $\tilde{\alpha}_{ni}$ are Adams-Bashforth coefficients corresponding to
the sequence of times~$\tU_n$ and $m^s(n)$ is defined by
$t^s_{m^s(n)} \le \tU_n < t^s_{m^s(n)+1}$, i.e., it is the index of the last
evaluation on set~$s$ that can influence the step $\dtU_n$ (see
Figure~\ref{fig:indexing}).  The interpolation coefficients are given by
\begin{equation}
  I_{ni;\allq} =
  \prod_{s=1}^S
  \ell_{q^s}\!\!\left(\tU_{n-i}; t^s_{m^s(n)}, \ldots,
    t^s_{m^s(n)-(k-1)}\right).
\end{equation}
For computational purposes, it is useful to rewrite these steps as
\begin{equation}\label{eq:lts-sum}
  \dty_n =
  \sum_{q^1=m^1(n)-(k-1)}^{m^1(n)} \cdots \sum_{q^S=m^S(n)-(k-1)}^{m^S(n)}
  \beta_{n;\allq}
  \D\ofyq,
\end{equation}
where the coefficient is
\begin{equation}\label{eq:lts-coef-small}
  \beta_{n;\allq} =
  \dtU_n \sum_{i=0}^{k-1} \tilde{\alpha}_{ni}
  \prod_{s=1}^S
  \ell_{m^s(n)-q^s}\!\!\left(\tU_{n-i};
    t^s_{m^s(n)}, \ldots, t^s_{m^s(n)-(k-1)}\right).
\end{equation}
The full change in the value of a given set of degrees of freedom over an
entire step can then be obtained by summing the contributions of all these
small steps.  This will give for each set of degrees of freedom an equation of
the form
\begin{equation}\label{eq:lts-form}
  \dy^s_m = \dt^s_m \sum_{q^1} \cdots \sum_{q^S} a^s_{m;\allq} \D^s\ofyq
\end{equation}
for some coefficients $a^s_{m;\allq}$.
Setting~\eqref{eq:lts-form} equal to
this sum over~\eqref{eq:lts-sum} gives the coefficients
\begin{equation}\label{eq:lts-coef-large}
  \dt^s_m a^s_{m;\allq} = \sum_{n=n^s(m)}^{n^s(m+1)-1} \beta_{n;\allq}.
\end{equation}

\begin{figure}
  \centering%
  \newcommand\stepline{\raisebox{0.5ex}{\rule{2em}{0.4pt}}}
  \begin{tabular}{
      M@{\hspace{1em}}M@{\hspace{1em}}McMcM@{\hspace{1em}}M@{\hspace{1em}}M}
    t^A \rightarrow \tU & \tU \rightarrow t^A
    &&&&&&
    \tU \rightarrow t^B & t^B \rightarrow \tU \\
    n^A(3) = 4 & m^A(4) = 3 &
    t^A_3 & \stepline & \tU_4 & \stepline & t^B_2 &
    m^B(4) = 2 & n^B(2) = 4 \\
    n^A(2) = 3 & m^A(3) = 2 &
    t^A_2 & \stepline & \tU_3 & & &
    m^B(3) = 1 & \\
    n^A(1) = 2 & m^A(2) = 1 &
    t^A_1 & \stepline & \tU_2 & & &
    m^B(2) = 1 & \\
          & m^A(1) = 0 &
    &     &             \tU_1 & \stepline & t^B_1 &
    m^B(1) = 1 & n^B(1) = 1 \\
    n^A(0) = 0 & m^A(0) = 0 &
    t^A_0 & \stepline & \tU_0 & \stepline & t^B_0 &
    m^B(0) = 0 & n^B(0) = 0 \\
  \end{tabular}
  \caption{Example of the values of $m^s(n)$ and $n^s(m)$ for an arbitrarily
    chosen step pattern on two sets.  These quantities give a mapping between
    the indices of the sequences of times $t^s_m$ and $\tU_n$, with $n^s(m)$
    mapping indices of $t^s_m$ to the corresponding indices of $\tU_n$ and
    $m^s(n)$ performing the reverse map.  In cases where there is no $t^s_m$
    corresponding to a given $\tU_n$ the index given by $m^s(n)$ is for the most
    recent step.}
  \label{fig:indexing}
\end{figure}

\section{Application to nearest-neighbor couplings}
\label{sec:special-cases}

\subsection{Element splitting}
\label{sec:element-splitting}

The equations~\eqref{eq:lts-form}
involve many more evaluations of the derivative than the
standard GTS Adams-Bashforth method, so in this form the LTS
method is unlikely to be more efficient.  However, if the couplings between the
sets of degrees of freedom are inexpensive to calculate compared to the
interactions within each set, then the required number of evaluations can be
reduced..  Let us suppose that the derivative on set~$s$ is split into a
``volume'' portion
depending only
on set~$s$ itself and a ``boundary'' portion
encoding the coupling to other sets:
\begin{equation}\label{eq:element-splitting}
  \D^s\ofyq = \V^s(\y^s_{q^s}) + \B^s\ofyq.
\end{equation}
These names are motivated by finite volume and discontinuous Galerkin methods,
where the terms from the interior and boundaries of elements split in this
manner.  Substituting \eqref{eq:element-splitting} into \eqref{eq:lts-sum} and
summing over the small
steps, the volume contribution to the full step on set~$s$ is
\begin{equation}
  (\dy^s_m)_{\text{vol}} =
  \mspace{-21mu}
  \sum_{q^s=m-(k-1)}^m
  \left[
    \sum_{n=n^s(m)}^{n^s(m+1)-1}
    \sum_{q^1=m^1(n)-(k-1)}^{m^1(n)}
    \mspace{-22mu} \cdots \mspace{5mu}
    \cancel{\sum_{q^s}}
    \mspace{5mu} \cdots \mspace{-22mu}
    \sum_{q^S=m^S(n)-(k-1)}^{m^S(n)}
    \mspace{-20mu}
    \beta_{n;\allq}
    \right]
  \V^s\!(\y^s_{q^s}),
\end{equation}
where $n^s(m)$ is defined by $\tU_{n^s(m)} = t^s_m$ (see
Figure~\ref{fig:indexing}).  This is the same form as the GTS Adams-Bashforth
method~\eqref{eq:gts} using the bracketed expression as coefficients (absorbing
the $\dt$ factor).  The bracketed expression does not depend on the form of the
derivative, so to evaluate it we can take the boundary coupling~$\B^s$ to be
zero, in which case this is the only contribution to the step.  As this is then
a $k$th-order GTS method and the Adams-Bashforth method is the unique
$k$th-order method of this form, the bracketed quantity must be the standard
Adams-Bashforth coefficient.  Returning to the general case with a coupling,
this shows that a set of degrees of freedom can be evolved using the standard
Adams-Bashforth method for the volume portion with only the coupling terms
evaluated using~\eqref{eq:lts-sum}.

This simplification applies in intermediate cases as well: if the full
derivative can be split into portions each of which depends on only some of the
degree-of-freedom sets, each of those contributions to the step can be
calculated independently using~\eqref{eq:lts-sum} ignoring non-contributing
sets.  In calculations where the sets are only coupled pairwise, this implies
that only the $S=2$ case need be considered.

\subsection{Two-set case}

As mentioned previously, a major intended application of this time-stepping
scheme is to evolution of PDEs using discontinuous Galerkin methods.
In such
a method, the sets of degrees of freedom are only coupled pairwise and
the update method reduces to a collection of standard Adams-Bashforth
methods and LTS methods with $S=2$.  For the two-set case, we call the sets~$A$
and $B$ and define the selection functions $\Theta^A_n$, $\Theta^B_n$, and
$\Theta^{AB}_n$ to be one if $\tU_n$ is an evaluation time for only set~$A$,
only set~$B$, or both sets, respectively.  By construction, a time evaluated on
neither set can never occur.  These selection functions sum to one, so we can
write $\beta_{n;q^Aq^B} = \beta^A_{n;q^Aq^B} + \beta^B_{n;q^Aq^B} +
\beta^{AB}_{n;q^Aq^B}$ with
\begin{multline}
  \beta^{A,B,AB}_{n;q^Aq^B} =
  \dtU_n \sum_{i=0}^{k-1} \tilde{\alpha}_{ni}
  \Theta^{A,B,AB}_{n-i}
  \ell_{m^A(n)-q^A}\!\!\left(\tU_{n-i};
    t^A_{m^A(n)}, \ldots, t^A_{m^A(n)-(k-1)}\right)
  \\\times
  \ell_{m^B(n)-q^B}\!\!\left(\tU_{n-i};
    t^B_{m^B(n)}, \ldots, t^B_{m^B(n)-(k-1)}\right).
\end{multline}
By the definition of $m(n)$, $\tU_n$ is not older than $t^{A,B}_{m^{A,B}(n)}$,
so, from the construction of the $\tU_n$ we see that $\tU_{n-i} \ge
t^{A,B}_{m^{A,B}(n) - (k-1)}$.  This implies that if $\tU_{n-i}$ is an
evaluation time for either set, it is one of the control points in the
corresponding Lagrange polynomial.  We can therefore collapse those polynomials
to obtain
\begin{align}
  \beta^A_{n;q^Aq^B} &=
  \Theta^A_{n^A(q^A)}
  \dtU_n \tilde{\alpha}_{n,n-n^A(q^A)}
  \ell_{m^B(n)-q^B}\!\!\left(t^A_{q^A};
    t^B_{m^B(n)}, \ldots, t^B_{m^B(n)-(k-1)}\right)
  \\&&\llap{$n - n^A(q^A) < k$} \notag
  \\
  \beta^{AB}_{n;q^Aq^B} &=
  \dtU_n \tilde{\alpha}_{n,n-n^A(q^A)}
  \qquad n^A(q^A) = n^B(q^B) \qquad n - n^A(q^A) < k,
\end{align}
and $\beta^B_{n;q^Aq^B} = \beta^A_{n;q^Bq^A}$.
The expressions should be taken to be zero when the conditionals on the right
are not satisfied.
Some example values are
shown in Tables~\ref{table:example-2to1} and \ref{table:example-2to1-start}.
The meaning of, for example, the first entry for (a) in
Table~\ref{table:example-2to1} is that in~\eqref{eq:lts-form} the coefficient
$a^A_{0;0,1} = 115/64$ (where we have chosen to number the steps starting from
$t=0$ so the step on set $B$ at $\dt^B$ is step 1), so the equation for this
step begins
\begin{equation}
  \y^A\bof{\dt^A} - \y^A\bof{0} =
  \dt^A \left[\frac{115}{64} \D^A\bof{\y^A\bof{0}, \y^B\bof{\dt^B}}
    + \cdots \right].
\end{equation}
Similarly, the lower-left entry for (c) in Table~\ref{table:example-2to1-start}
indicates that one term in the second small step is
\begin{equation}
  \y^B\bof{2\dt^B} - \y^B\bof{\dt^B} =
  \dt^B \left[\frac{2}{3} \D^A\bof{\y^A\bof{-2\dt^A}, \y^B\bof{\dt^B}}
    + \cdots \right].
\end{equation}
Additional tables of coefficients can be found in Appendix~\ref{sec:tables}.

\begin{table}
  \centering%
  \raisebox{-0.5\height}{%
    \begin{tikzpicture}
      \Astep{-2}{-2\dt^A}
      \Astep{-1}{-\dt^A}
      \Astep{0}{0}
      \Astep{1}{\dt^A}
      \Bstepunused{-4/2}
      \Bstepunused{-3/2}
      \Bstep{-2/2}{-2\dt^B}
      \Bstep{-1/2}{-\dt^B}
      \Bstep{0}{0}
      \Bstep{1/2}{\dt^B}
      \Bstep{2/2}{2\dt^B}
      \Achart{1/2}{(a)}
      \Bchart{1/4}{(b)}
      \Bchart{3/4}{(c)}
    \end{tikzpicture}%
  }%
  \hspace{1em}%
  \begin{tabular}{M|MMMM}
    \multicolumn{1}{c|}{\text{(a)}}
            & \dt^B           & 0            & -\dt^B         & -2\dt^B      \\
    \hline
    0       & \frac{115}{64}  & \frac{7}{24} & -\frac{11}{64} & 0            \\
    -\dt^A  & -\frac{115}{96} & 0            & -\frac{11}{32} & \frac{5}{24} \\
    -2\dt^A & \frac{23}{64}   & 0            & \frac{11}{192} & 0            \\
    \multicolumn{5}{l}{} \\
    \multicolumn{1}{c|}{\text{(b)}}
            & \dt^B & 0             & -\dt^B       & -2\dt^B      \\
    \hline
    0       &       & \frac{23}{12} & -\frac{1}{2} & 0            \\
    -\dt^A  &       & 0             & -1           & \frac{5}{12} \\
    -2\dt^A &       & 0             & \frac{1}{6}  & 0            \\
    \multicolumn{5}{l}{} \\
    \multicolumn{1}{c|}{\text{(c)}}
            & \dt^B           & 0            & -\dt^B        & -2\dt^B \\
    \hline
    0       & \frac{115}{32}  & -\frac{4}{3} & \frac{5}{32}  &         \\
    -\dt^A  & -\frac{115}{48} & 0            & \frac{5}{16}  &         \\
    -2\dt^A & \frac{23}{32}   & 0            & -\frac{5}{96} &         \\
  \end{tabular}
  \caption{A third-order method for two sets~$A$ and $B$ with $B$ evaluated
    twice as often as $A$.  Coefficients $a^s_{m;q^Aq^B}$
    in~\eqref{eq:lts-form} for use with derivatives
    evaluated using data from $A$ and $B$ at the times
    indicated for (a) a step of set~$A$ from $0$ to $\dt^A$, and steps of
    set~$B$ (b) from $0$ to $\dt^B$ and (c) from $\dt^B$ to $2\dt^B = \dt^A$.}
  \label{table:example-2to1}
\end{table}

\begin{table}
  \centering%
  \raisebox{-0.5\height}{%
    \begin{tikzpicture}
      \Astep{-2}{-2\dt^A}
      \Astep{-1}{-\dt^A}
      \Astep{0}{0}
      \Astep{1}{\dt^A}
      \Astepunused{2}
      \Bstep{-2}{-2\dt^A}
      \Bstep{-1}{-\dt^A}
      \Bstep{0}{0}
      \Bstep{1/2}{\dt^B}
      \Bstep{2/2}{2\dt^B}
      \Bstepunused{3/2}
      \Bstepunused{4/2}
      \Achart{1/2}{(a)}
      \Bchart{1/4}{(b)}
      \Bchart{3/4}{(c)}
    \end{tikzpicture}%
  }%
  \hspace{1em}%
  \begin{tabular}{M|MMMM}
    \multicolumn{1}{c|}{\text{(a)}}
            & \dt^B         & 0           & -\dt^A       & -2\dt^A      \\
    \hline
    0       & \frac{5}{3}   & \frac{1}{4} & 0            & 0            \\
    -\dt^A  & -\frac{10}{9} & 0           & -\frac{2}{9} & 0            \\
    -2\dt^A & \frac{1}{3}   & 0           & 0            & \frac{1}{12} \\
    \multicolumn{5}{l}{} \\
    \multicolumn{1}{c|}{\text{(b)}}
            & \dt^B & 0             & -\dt^A        & -2\dt^A     \\
    \hline
    0       &       & \frac{17}{12} & 0             & 0           \\
    -\dt^A  &       & 0             & -\frac{7}{12} & 0           \\
    -2\dt^A &       & 0             & 0             & \frac{1}{6} \\
    \multicolumn{5}{l}{} \\
    \multicolumn{1}{c|}{\text{(c)}}
            & \dt^B         & 0              & -\dt^A       & -2\dt^A \\
    \hline
    0       & \frac{10}{3}  & -\frac{11}{12} & 0            &         \\
    -\dt^A  & -\frac{20}{9} & 0              & \frac{5}{36} &         \\
    -2\dt^A & \frac{2}{3}   & 0              & 0            &         \\
  \end{tabular}
  \caption{Rules for reducing the time step size in one set to start the
    algorithm in Table~\ref{table:example-2to1} from a GTS state.  For $t \le
    0$ both sets step together at interval~$\dt^A$, after which set~$B$ changes
    to a step of $\dt^B = \dt^A/2$.  For steps before $t = 0$ the standard GTS
    rules can be used, and for steps beyond $t = \dt^A$ the rules in
    Table~\ref{table:example-2to1} apply.}
  \label{table:example-2to1-start}
\end{table}

\section{Numerical results}
\label{sec:numerical}

We tested this scheme on a set of field equations evaluated using discontinuous
Galerkin methods.  In a DG formulation, the domain of evolution is divided into
elements, with each element containing a collection of nodes.  The evolution
equations are evaluated locally within each element and this collection of
partial solutions is coupled by adding additional terms at the element
boundaries obtained from comparison with neighboring elements.
The problem therefore naturally splits as described in
Section~\ref{sec:element-splitting}, with each element being evolved using a
standard GTS method and the couplings using the LTS equations.

\subsection{Convergence}
\label{sec:numerical:convergence}

For our first test problem, we wanted a simple nonlinear
problem to test the
convergence and conservation properties of the method.  Accordingly,
we evolved the Burgers equation:
\begin{equation}\label{eq:burgers}
  0 = \pderiv[u]{t} + \pderiv{x}\left(\frac{1}{2} u^2\right).
\end{equation}
The DG elements were coupled using the numerical flux of Harten, Lax, and van
Leer~\cite{Harten1983}.  In the common case where
the solution has the same sign on
both sides of the interface
between two elements, this reduces to an upwind
flux.  For
this PDE, the spatial integral of the field~$u$ is a linear conserved quantity,
and
this carries over to the discretized system when using a DG scheme.

For our first test, we considered the exact
solution
\begin{equation}\label{eq:bump}
  u(t ,x)
  = \frac{\sqrt{1 - 4 t (x - t)} - 1 + 2 t x}{2 t^2}
  = \frac{
    2 \big(\sqrt{1 - 4 t (x - t)} + 1 - 2 x (x - t)\big)
  }{
    \big(\sqrt{1 - 4 t (x - t)} + 1\big)^2
  }.
\end{equation}
We used the analytic value of $u(t=-1/8, x)$ as an initial
condition and evolved $u$
from $t = -1/8$ to $t = 3/2$ on the interval $-9/8 \le x \le 1/8$,
with free boundary conditions.  The
interval was divided into 16 elements of equal size, each of which contained 10
nodes distributed as Legendre-Gauss-Lobatto points.
These values reduce the spatial truncation error to approximately
$10^{-15}$, smaller than the floating-point roundoff error, so it can be
ignored for this analysis.

\begin{figure}
  \centering%
  \input{plots/step_pattern_burgers}
  \caption{
    Stepping pattern used to evolve \eqref{eq:bump} using a 5th-order
    integrator using the step-size condition $u \dt < 2^{-12}$.  The small
    steps used at the start of the integration are not visible at this scale.}
  \label{figure:stepping-bump}
\end{figure}
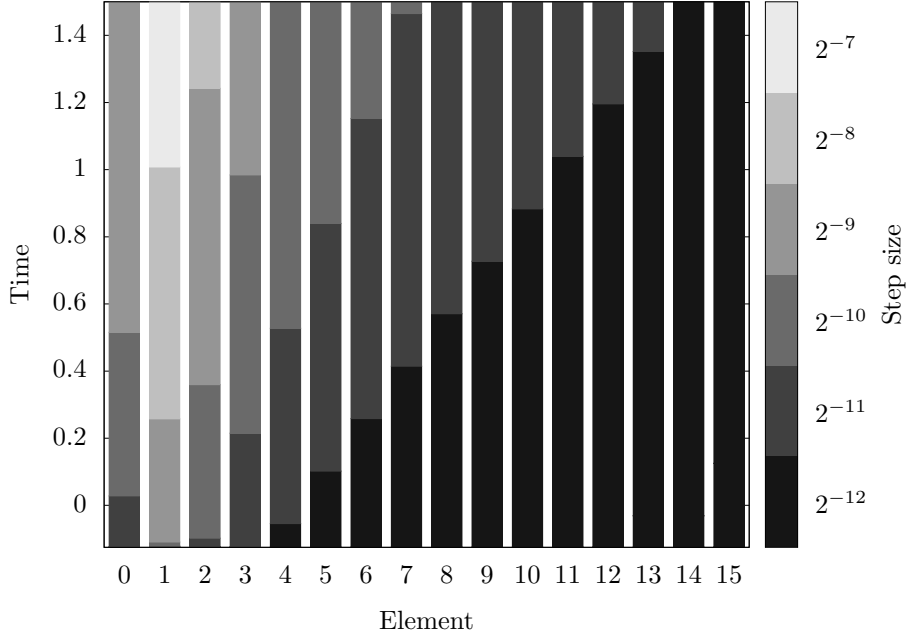

To simplify bookkeeping,
the size of each time step was restricted to be a power of two, but each
element was allowed to independently adjust its step size to conform to the CFL
limit.  We found that increasing the step size too rapidly led to a growth in
the error, so each step size increase was limited to a factor of two and an
increase was only allowed if the previous $k-1$ steps were all of the same
size, where $k$ is the integration order.  All elements were initialized with a
very small time step of $2^{-27}$.  Figure~\ref{figure:stepping-bump} shows a
typical step pattern.

\begin{figure}
  \centering%
  \input{plots/convergence}
  \caption{
    Convergence of the evolution of \eqref{eq:bump} with decreasing step
    size for integrators of the indicated expected orders.  The slopes of the
    straight lines indicate the expected rates of convergence.  The domain of
    evolution was divided into 16 elements, each of which has an independently
    chosen step size based on the local CFL condition, which changed
    dynamically throughout the evolution.  A typical pattern of steps is shown
    in Figure~\ref{figure:stepping-bump}.  The visible oscillations are due to
    our restriction that the time steps must always be powers of two: a
    uniformly scaled-down step pattern occurs only after the step sizes have
    been reduced by half.}
  \label{figure:convergence}
\end{figure}
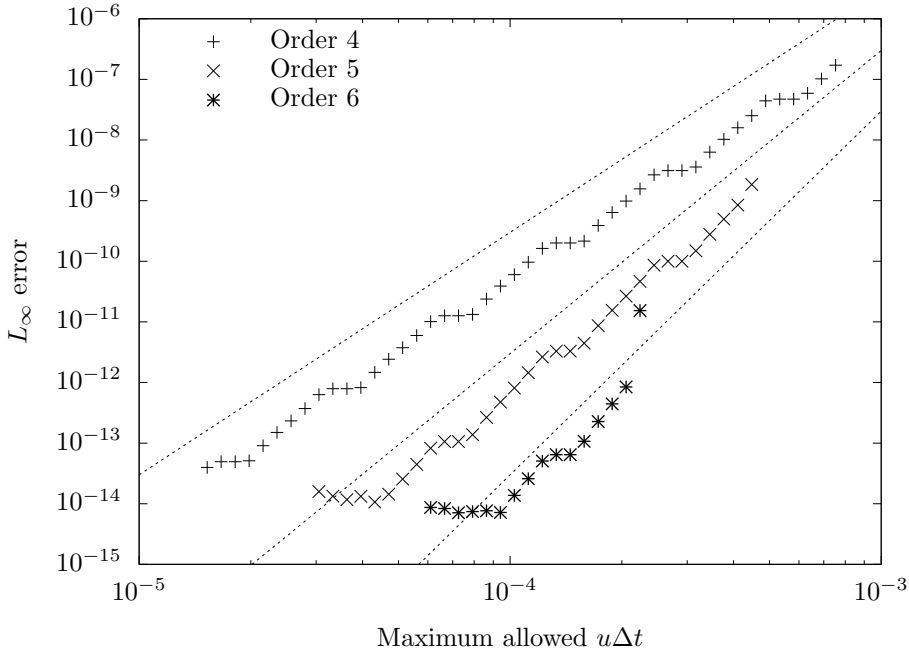

The sole characteristic speed in the Burgers system is $\abs{u}$, so the CFL
condition predicts that the stable step size will be proportional to $1/u$.
The coefficient of proportionality can be estimated from the point
distribution and properties of the time stepper, but in order to investigate
the convergence of the method we kept it as a free parameter.  The convergence
as a function of this parameter is shown in Figure~\ref{figure:convergence}.

\subsection{Conservation}
\label{sec:numerical:conservation}

\begin{figure}
  \centering%
  \input{plots/wave}
  \caption{
    Numerical evolution of the periodic solution with initial condition
    given by \eqref{eq:wave}.  The exact solution forms a shock at $t \approx
    0.37$, after which time the numerical solution becomes qualitatively
    incorrect.}
  \label{figure:wave-evolution}
\end{figure}
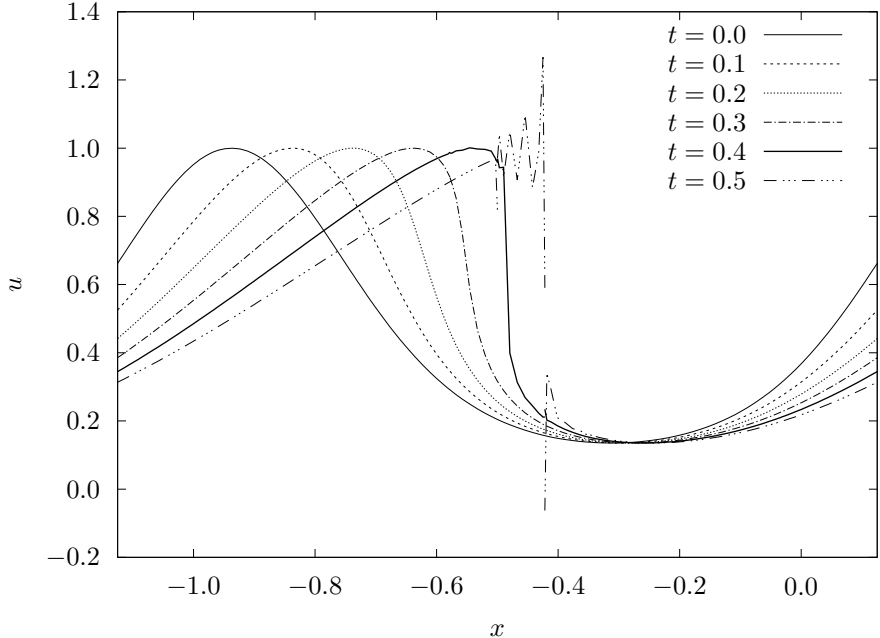

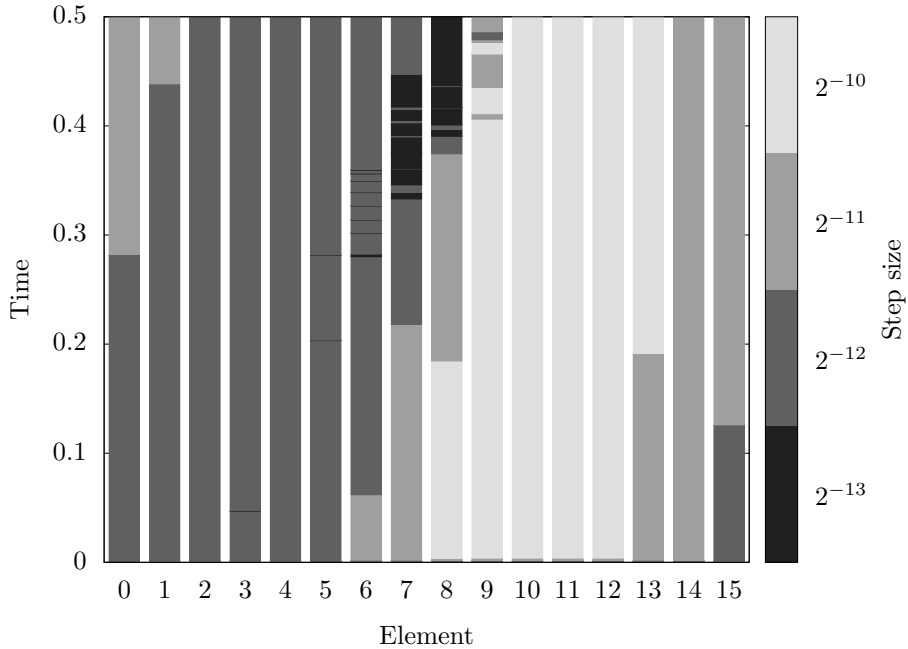
\begin{figure}
  \centering%
  \input{plots/step_pattern_periodic}
  \caption{
    Stepping pattern used to evolve \eqref{eq:wave} using a 5th-order
    integrator using the step-size condition $u \dt < 2^{-12}$.  The small
    steps used at the start of the integration are not visible at this scale.
    The fluctuating step sizes at later times are due to the inability of the
    spatial method to resolve the solution as a shock forms at $t \approx
    0.37$.}
  \label{figure:stepping-wave}
\end{figure}

\begin{figure}
  \centering%
  \input{plots/conservation}
  \caption{
    Change in the integral of the conserved quantity~$u$ over the domain for
    the evolution of the wave given in~\eqref{eq:wave}.  Note that the
    timescale shown here is much larger than the time required for an
    (unresolved) shock to form.}
  \label{figure:conservation}
\end{figure}
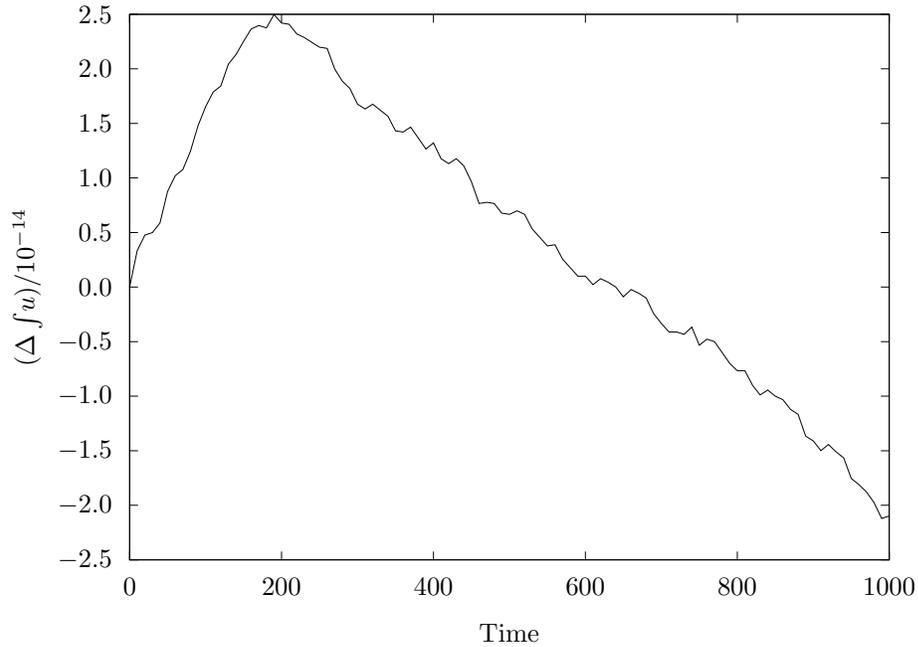

The previous test does not demonstrate conservation because of flow through the
domain boundaries.  To eliminate this effect, we perform an evolution of the
Burgers equation on a periodic domain.  We use the same domain decomposition as
for the previous example, but identify the endpoints at $-9/8$ and $1/8$.  We
initialize the field at $t=0$ to
\begin{equation}\label{eq:wave}
  u(x) = \frac{1}{e} \exp\left(\sin\left(\frac{8 \pi}{5} x\right)\right).
\end{equation}
The exact evolution of these initial conditions forms a shock at
\begin{equation}
  t_s = \frac{5 e}{8 \pi} \frac{1}{
    \exp\left(\frac{\sqrt{5} - 1}{2}\right) \sqrt{\frac{\sqrt{5} - 1}{2}}}
  \approx 0.37,
\end{equation}
and the numerical evolution becomes poorly behaved at approximately that time
as shown in Figure~\ref{figure:wave-evolution}.  This results in the
fluctuating time steps visible in Figure~\ref{figure:stepping-wave}.  Despite
the large departure from the correct solution, we find that the integral of $u$
is conserved to
within a few orders of magnitude of machine epsilon, as shown in
Figure~\ref{figure:conservation}.

\subsection{Speed}
\label{sec:numerical:speed}

\begin{figure}
  \centering%
  \begin{tikzpicture}
    \draw (-3.5, -3.5) grid (4.5, 4.5);
    \newcommand\refined[3]{
      \pgfmathsetmacro\numlines{2^#1-1}
      \foreach \x in {1, ..., \numlines} {
        \pgfmathsetmacro\blockfrac{\x / 2^#1}
        \draw (#2, #3 + \blockfrac) -- (#2 + 1, #3 + \blockfrac);
        \draw (#2 + \blockfrac, #3) -- (#2 + \blockfrac, #3 + 1);
      }
    }
    \foreach \r in {-3, ..., 3} {
      \foreach \c in {-3, ..., 3} {
        \pgfmathtruncatemacro\refinement{4 - (abs(\r) + abs(\c))}
        \ifnum\refinement>0{
          \refined{\refinement}{\r}{\c}
        }\fi
      }
    }
  \end{tikzpicture}
  \caption{
    Cross-section through the center of the evolution domain used for
    the 3D~gauge wave speed test in Section~\ref{sec:numerical:speed}.  There
    are 31 additional rings of elements the same size as the largest elements
    shown here surrounding this $7\times7$ section.  In three dimensions, the
    refined region in the center has an octahedral shape.}
  \label{figure:gauge-wave-domain}
\end{figure}
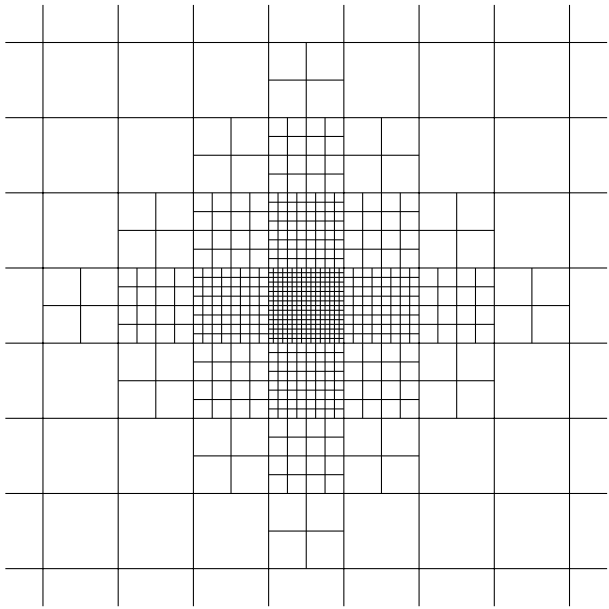

To measure the speed improvement achieved by using the LTS method, we examined
a system of more practical interest than the Burgers equation: the generalized
harmonic formulation of the Einstein field equations~\cite{Lindblom_2006} with
elements coupled by the upwind flux given by (6.3) of~\cite{Porth2017}.  We
evolve the gauge wave solution given in Section~4.2 of \cite{Alcubierre_2003},
deriving the evolved quantities from the metric
\begin{gather}
  \psi_{ab} = \operatorname{diag}(-H, H, 1, 1) \label{eq:gauge-wave} \\
  H = 1 - A \sin\left(\frac{2 \pi (x - t)}{d}\right)
\end{gather}
with amplitude~$A = 0.01$ and wavelength~$d = 0.125$.  The evolution domain
consisted of a periodic unit cube divided into $69^3$ equal-sized blocks.  The
central block was refined into $16^3$ equal-sized elements, and nearby blocks
were refined uniformly as needed to keep the lengths of adjacent elements
within a factor of two of each other, as shown in
Figure~\ref{figure:gauge-wave-domain}.

The CFL step size bound on the largest elements is 16~times the bound on the
smallest elements, so if the cost of the simulation is dominated by the large
elements, the LTS algorithm could achieve a theoretical speedup of 16.  For the
domain size chosen, however, there are sufficiently many small elements that
cannot be improved by the LTS algorithm that the theoretical speedup is reduced
to approximately $12.7$.

These evolutions were performed on a cluster of 126 Lenovo
Nextscale NX360 M5 compute nodes, each with 24 cores and 64~GB of RAM.  The
nodes communicated using Mellanox Infiniband.

We evolved the wave using integrators with orders from 1 to 8.  For
orders~1--4, evolution was performed on 8 nodes.
For orders~5--8, the simulation was run on 12 nodes because of the higher
memory requirements.  The first order integrator was run to time~0.28,
resulting in a total run time for the GTS evolution of approximately 10~hours.
The final times for the higher-order evolutions were adjusted to keep the run
times of the GTS evolutions similar by taking into account the Adams-Bashforth
stability factors (Table~\ref{table:ab-stability}) and the number of nodes
used.

We found that the LTS evolutions of the wave ran 6--9 times faster than the
corresponding GTS runs.  These results were not strongly dependent on the order
of the method.

\subsection{Stability}

\begin{table}
  \centering%
  \newcommand\colwidth{2.3em}%
  \newcolumntype{z}{>{\centering\arraybackslash$}m{\colwidth}<{$}}
  \begin{tabular}{r|zzzzzzzz}
    Order & 1 & 2 & 3 & 4 & 5 & 6 & 7 & 8 \\
    \hline
    $C$ &
    1 &
    \frac{1}{2} &
    \frac{3}{11} &
    \frac{3}{20} &
    \frac{45}{551} &
    \frac{5}{114} &
    \frac{945}{40663} &
    \frac{945}{77432}
  \end{tabular}
  \caption{Stability factors $C$ in \eqref{eq:stability-factor} for
    advection problems using Adams-Bashforth methods on finite-difference
    grids.}
  \label{table:ab-stability}
\end{table}

For common GTS integrators, integration of a PDE will be stable if the step
size taken is bounded by
\begin{equation}\label{eq:stability-factor}
  \dt < C \frac{\dx}{\charspeed},
\end{equation}
where $\dx$ is the minimum spacing in the spatial discretization, $\charspeed$
is the largest characteristic speed of the PDE, and $C$ is a constant that
depends on the system being integrated and the choice of integrator.  For the
standard Adams-Bashforth integrators, the values of $C$ for solving advection
problems with finite-difference methods are given in
Table~\ref{table:ab-stability}.  The values of $C$ for other systems are
usually somewhat smaller, so we define
\begin{equation}
  s = \frac{C}{C_{\text{advection}}}.
\end{equation}

We measured the stability of our integrators when evolving the one-dimensional
generalized harmonic equations using the gauge wave
solution~\eqref{eq:gauge-wave} with the transverse dimensions suppressed.  We
chose a wavelength of $d=2\pi$ and an amplitude of $A=10^{-6}$.  A small
amplitude was chosen so the wave would have a negligible effect on the
characteristic speeds of the system while still ensuring any unstable modes
were excited.

We evolved the wave on a periodic domain of length~$4\pi$.  This was divided
into two equal segments, one of which was refined into~$2^m$ equal elements,
and the other of which was refined into~$2^n$ with $n \ge m \ge 2$.  The
solution on each element was described by 3~Legendre-Gauss-Lobatto points.

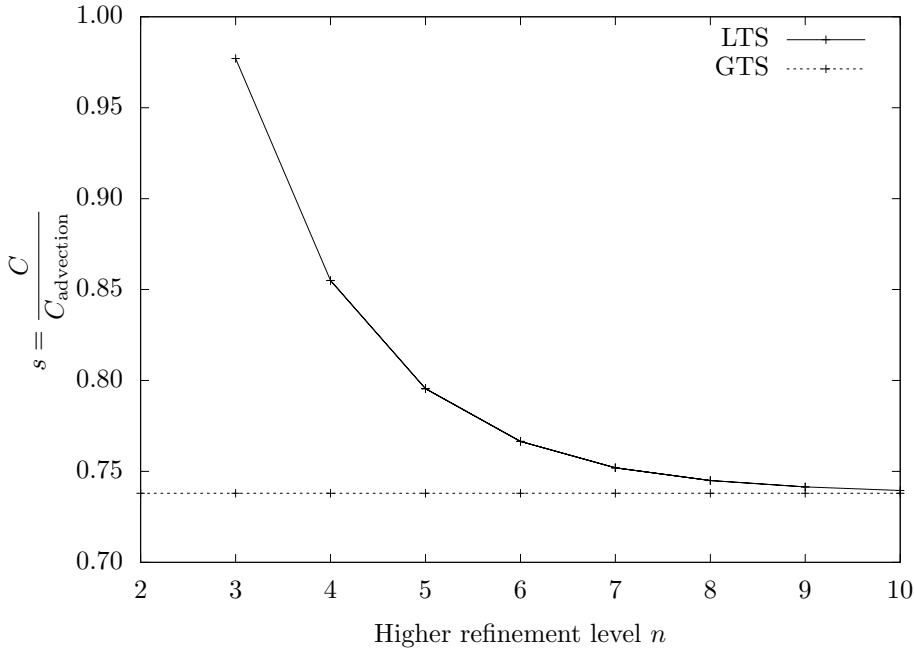
\begin{figure}
  \centering%
  \input{plots/lts_cfl}
  \caption{
    Stable time step for integration of the 1D~gauge wave, given as the CFL
    factor $C$ in \eqref{eq:stability-factor} divided by the ``advection''
    values in Table~\ref{table:ab-stability}.  The domain is divided into two
    parts with $2^m$ and $2^n$ elements with $n \ge m \ge 2$.  For LTS
    integration, element size differences of up to $n - m = 3$ are shown, but
    the curves are indistinguishable at our measurement tolerance.  These
    values were observed to be independent of integrator order for order at
    least 2.}
  \label{figure:lts-cfl}
\end{figure}

We found the first-order time stepper to be very unstable ($s < 0.1$) for both
global and local time-stepping.  For steppers of order~2 and above, we found
the values of $s$ to be independent of time-stepper order to within our
measurement error of approximately $0.001$.  When run with global time-stepping
($m = n$), we measured $s = 0.738$, independent of the refinement level.  With
non-uniform refinement ($n > m$) the LTS integrator was found to be more stable
than the GTS integrator, with $s$ approaching the GTS value as the number of
elements was increased, as shown in Figure~\ref{figure:lts-cfl}.  For element
size ratios from $n - m = 1$ to $n - m = 3$, we found that, to within our
measurement error, $s$ was independent of the size of the large elements.

\subsection{Convergence at constant CFL factor}

Section~\ref{sec:numerical:convergence} considered the convergence of the
method with a particular discretization of the computational domain.  When
solving PDEs, another relevant test is the convergence as the domain is
spatially refined and the time step is adjusted as necessary to maintain
stability.  We measured this rate using the one-dimensional gauge wave problem
described in the previous section with a $2:1$ ratio of element sizes ($n = m +
1$).  The time steps were restricted to be powers of 2 and their values were
chosen using the CFL bound with $s = 0.5$.  We found that, in most cases, the
dominant source of error was the spatial discretization, and we were only able
to make accurate measurements of the error caused by the time stepper for
low-order integrators on comparatively high-order elements.

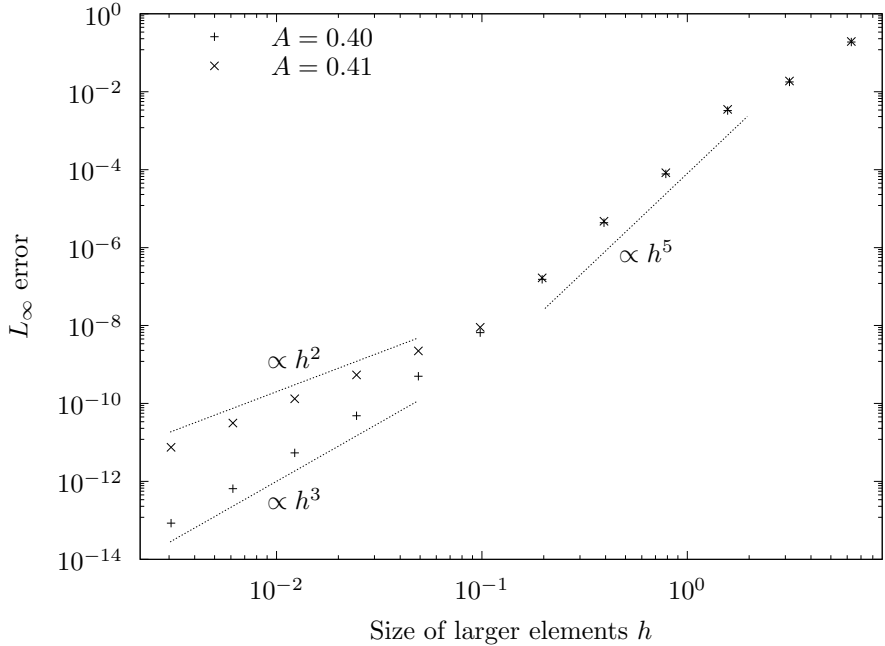
\begin{figure}
  \centering%
  \input{plots/resolution_convergence}
  \caption{
    Error in the value of the one-dimensional gauge-wave~\eqref{eq:gauge-wave}
    at $t=1$ with $2:1$ spatial refinement of elements, evolved on fourth-order
    elements with a third-order time stepper.  The step size was chosen to be a
    power of two and required to satisfy \eqref{eq:stability-factor} with
    $s=0.5$.  This results in a static step pattern for wave amplitudes below
    $A\approx0.4008$ and a time-varying step pattern for larger amplitudes.}
  \label{figure:cfl-convergence}
\end{figure}

The largest characteristic speed of the gauge wave
solution~\eqref{eq:gauge-wave} is $\sqrt{H}$.  For small amplitude waves, this
speed will be nearly constant, but for larger amplitudes there can be
significant variation as the wave propagates through the domain.  At some
threshold amplitude, this will start to cause time-dependent variation in the
step size by causing the limit from \eqref{eq:stability-factor} to cross one of
the permitted power-of-two step sizes.  We found that, when the amplitude was
small enough that the step pattern was static, the convergence rate was
unaffected by the step size variation across the domain.  When the amplitude
was large enough to cause a time-varying step pattern, however, the convergence
order was reduced by one.  The convergence for the third-order integrator is
shown in Figure~\ref{figure:cfl-convergence}.

\section{Conclusions}

When integrating systems of coupled ODEs, particularly those arising from
discretizations of PDE systems, it is often the case that time-step-related
restrictions arise primarily in a small subset of the variables being
integrated.  Using standard evolution schemes, this forces all degrees of
freedom to be evolved with the most restrictive stable time step, potentially
causing significant waste of computational resources.  A local time-stepping
integrator removes this requirement, allowing different degrees of freedom to
be updated at different frequencies.

This paper has presented a local time-stepping scheme based on the
Adams-Bashforth family of multistep integrators.  This method allows arbitrary
step choices, with a completely independent choice of time step for each
variable.  Unlike some previous schemes, it retains the full convergence order
of the Adams-Bashforth integrator it is based on.  This method is also
conservative in that all linear conserved quantities of the system are constant
to numerical roundoff under evolution.
The method will be efficient for DG methods or other algorithms with only local
spatial couplings.

The use of this method was demonstrated on
several systems evolved using
the DG framework.  The roundoff-level conservation of a conserved quantity was
observed, and the expected convergence rate
as the time step was decreased was
observed for multiple integrator orders.
The method was observed to converge one order more slowly when applied
to a system with decreasing mesh spacing and dynamically updated time step
sizes.
For the evolution of a gauge wave using the
generalized harmonic equations, we observed an
evolution speed improvement by
up to a factor of 6--9
from switching
from the global to the local scheme.

This method will be used for DG evolutions of general relativity and
magnetohydrodynamics in upcoming work using the SpECTRE
code~\cite{SpECTRE2017}.

\appendix

\section{Consistency and order of convergence}
\label{sec:consistency}

\newcommand\f{\mat{f}}
\newcommand\order[1]{O(\Delta t^{#1})}
\newcommand\jlimit{\substack{j^1,\ldots,j^S=0\\j^1+\cdots+j^S<k}}

We will demonstrate here that the method given by \eqref{eq:lts-sum} and
\eqref{eq:lts-coef-small} is consistent at the same order~($k$) as its base
method.
If we assume that all past values have been determined accurately, with error
at most $\order{k}$, then we can write
\begin{equation}
  \y^s_q = \y^s(t^s_q) + \order{k}.
\end{equation}
If we define a function
\begin{equation}
  \f\bof{t^1, \ldots, t^S} = \D\bof{\y^1\bof{t^1}, \ldots, \y^S\bof{t^S}},
\end{equation}
then, assuming $\D$ is sufficiently smooth ($C^{k-1}$), we can expand
$\D$ for steps near
time $\tU_n$ as
\begin{align}
  \D\ofyq &= \f\bof{t^1_{q^1}, \ldots, t^S_{q^S}} + \order{k} \\
  &= \sum_{\jlimit}^{k-1} \f_{j^1\cdots j^S}
  \frac{{(t^1_{q^1} - \tU_n)}^{j^1} \cdots {(t^S_{q^S} - \tU_n)}^{j^S}}{
    j^1! \cdots j^S!}
  + \order{k},
  \label{eq:deriv-Taylor}
\end{align}
where $\f_{j^1\cdots j^S}$ is the value of $\f$ differentiated $j^1$ times with
respect to its first argument, $j^2$ times with respect to its second, etc.,
and evaluated at $(t^1, \ldots, t^S) = (\tU_n, \ldots, \tU_n)$.

If we now substitute \eqref{eq:lts-coef-small} and \eqref{eq:deriv-Taylor} into
\eqref{eq:lts-sum}, we find that the entire expression factors into separate
parts for each set of degrees of freedom:
\begin{multline}\label{eq:consistency:factored}
  \dty_n =
  \dtU_n \sum_{i=0}^{k-1} \tilde{\alpha}_{ni}
  \sum_{\jlimit}^{k-1} \f_{j^1\cdots j^S}
  \\ \times
  \prod_{s=1}^S \left[
    \sum_{q=0}^{k-1}
    \frac{{(t^s_{m^s(n) - q} - \tU_n)}^{j^s}}{j^s!}
    \ell_q\!\!\left(\tU_{n-i};
    t^s_{m^s(n)}, \ldots, t^s_{m^s(n)-(k-1)}\right)
    \right]
  + \order{k+1}.
\end{multline}
The bracketed expression can be interpreted as an order-$(k-1)$ Lagrange
interpolation from the times $t^s_{m^s(n)}, \ldots, t^s_{m^s(n)-(k-1)}$ to $t =
\tU_{n-i}$ of the polynomial $(t - \tU_n)^{j^s}/j^s!$.  As $j^s \le k-1$, this
interpolation is exact and we can simplify the above expression to
\begin{equation}
  \dty_n =
  \dtU_n \sum_{i=0}^{k-1} \tilde{\alpha}_{ni}
  \sum_{\jlimit}^{k-1} \f_{j^1\cdots j^S}
  \frac{(\tU_{n-i} - \tU_n)^{j^1 + \cdots + j^S}}{j^1! \cdots j^S!}
  + \order{k+1}.
\end{equation}
The sum over the $j^s$ can now again be interpreted as a Taylor series for
$\f$, but evaluated this time with all arguments $\tU_{n-i}$:
\begin{equation}
  \dty_n
  = \dtU_n \sum_{i=0}^{k-1} \tilde{\alpha}_{ni}
  \f\bof{\tU_{n-i}, \ldots, \tU_{n-i}}
  + \order{k+1}
  = \dtU_n \sum_{i=0}^{k-1} \tilde{\alpha}_{ni} \deriv[\y]{t}\bof{\tU_{n-i}}
  + \order{k+1},
\end{equation}
where the second equality follows from the definition of $\f$.  This is now the
expression for a single step of an Adams-Bashforth integral.  The step has
error $\order{k+1}$, so it is equal to its analytic value to the accuracy of
this expression.  We can therefore perform the integral to find
\begin{equation}
  \dty_n = \y\bof{\tU_{n+1}} - \y\bof{\tU_n} + \order{k+1}.
\end{equation}

\section{Error term}
\label{sec:error}

By repeating the calculation of consistency above keeping an additional term
(and assuming now that $\D$ has one more continuous derivative) we can obtain
an estimate of the error of the method.  If we take the old values to be
accurate to order~$k$ (which, for an order-$k$ method, we can assume them to
be)
then the calculations in the previous section are nearly unchanged until
we reach~\eqref{eq:consistency:factored}, which now reads
\begin{multline}
  \dty_n =
  \dtU_n \sum_{i=0}^{k-1} \tilde{\alpha}_{ni}
  \sum_{\substack{j^1,\ldots,j^S=0\\j^1+\cdots+j^S \le k}} \f_{j^1\cdots j^S}
  \\ \times
  \prod_{s=1}^S \left[
    \sum_{q=0}^{k-1}
    \frac{{(t^s_{m^s(n)-q} - \tU_n)}^{j^s}}{j^s!}
    \ell_q\!\!\left(\tU_{n-i};
    t^s_{m^s(n)}, \ldots, t^s_{m^s(n)-(k-1)}\right)
    \right]
  + \order{k+2}.
\end{multline}
As before, the bracketed quantity is an order-$(k-1)$ polynomial interpolation,
but now the interpolations with $j^s = k$ are not exact.  Evaluating the
interpolations (and noting that $j^s = k$ for some set implies that $j^s = 0$
for all the others) results in
\begin{multline}
  \dty_n =
  \dtU_n \sum_{i=0}^{k-1} \tilde{\alpha}_{ni}
  \Biggg\{
  \sum_{\substack{j^1,\ldots,j^S=0\\j^1+\cdots+j^S \le k}} \f_{j^1\cdots j^S}
  \frac{{(\tU_{n-i} - \tU_n)}^{j^1+\cdots+j^S}}{j^1!\cdots j^S!}
  \\
  +
  \sum_{s=1}^S \frac{\partial^k\f}{(\partial t^s)^k}
  \mspace{-5mu}
  \left[
    \sum_{q=0}^{k-1}
    \frac{{(t^s_{m^s(n)-q} - \tU_n)}^k}{k!}
    \ell_q\!\!\left(\tU_{n-i};
    t^s_{m^s(n)}, \ldots, t^s_{m^s(n)-(k-1)}\right)
    - \frac{{(\tU_{n-i} - \tU_n)}^k}{k!}
    \right]
  \mspace{-8mu}
  \Biggg\}
  \\
  + \order{k+2}.
\end{multline}
The first term in the braces is the Taylor expansion of $d\y/dt$ expanded in
terms of $\f$.  Using the sum over $i$, we can integrate it to get $\y$ and a
correction term because the Adams-Bashforth integration does not integrate the
order-$k$ term exactly.  In the second term, the sum over $i$ can be used to
integrate the Lagrange polynomial exactly.  These manipulations produce
\begin{multline}
  \dty_n =
  \y\bof{\tU_{n+1}} - \y\bof{\tU_n}
  +
  \derivn[\y]{t}{k+1} \left[
    \dtU_n \sum_{i=0}^{k-1} \tilde{\alpha}_{ni}
    \frac{{(\tU_{n-i} - \tU_n)}^k}{k!}
    -
    \frac{\dtU_n^{k+1}}{(k+1)!}
    \right]
  \\
  +
  \sum_{s=1}^S \frac{\partial^k\f}{(\partial t^s)^k}
  \left[
    \dtU_n \sum_{q=0}^{k-1} \bar{\alpha}^s_{nq}
    \frac{{(t^s_{m^s(n)-q} - \tU_n)}^k}{k!}
    -
    \dtU_n \sum_{i=0}^{k-1} \tilde{\alpha}_{ni}
    \frac{{(\tU_{n-i} - \tU_n)}^k}{k!}
    \right]
  + \order{k+2},
  \mspace{-3mu}
\end{multline}
where we have defined the quantities
\begin{equation}
  \bar{\alpha}^s_{nq} =
  \frac{1}{\dtU_n} \int_{\tU_n}^{\tU_{n+1}} dt\,
  \ell_q\!\!\left(t; t^s_{m^s(n)}, \ldots, t^s_{m^s(n)-(k-1)}\right),
\end{equation}
which are the Adams-Bashforth coefficients for integrating from $\tU_n$ to
$\tU_{n+1}$ using known derivatives at the evaluation times on set~$s$.

We see now that there are two contributions to the error.  The first has no
dependence on the split of the solution into sets of degrees of freedom, and is
therefore the error for the GTS Adams-Bashforth method.  The second term is
specific to the LTS method and occurs because of the difference in the error
from
using values from the different step patterns.

\section{Element splitting for general methods}
\label{sec:general-splitting}

When comparing integrators, one may wish to use a GTS integrator that is not
usually expressed in terms of volume terms and boundary couplings (for example,
a Runge-Kutta method) in a framework designed for an LTS integrator that is so
expressed.  This is easiest if the GTS integrator can be cast into the element
splitting form (Section~\ref{sec:element-splitting}).

All common explicit GTS integrators (both multistep and substep) can be written
in the form
\begin{equation}\label{eq:general-framework:gts}
  u_{n+1} - u_n = \sum_i A_n^i (u_n - u_{n-i}) + \dt_n \sum_i B_n^i D(u_{n-i}).
\end{equation}
Adams-Bashforth integrators are usually written in this form with the $A_n^i =
0$.  Runge-Kutta methods take some manipulation.  For example, the second-order
midpoint method
\begin{equation}
  u_{n+1} - u_n = \dt\, D\left(u_n + \frac{1}{2} \dt D(u_n)\right)
\end{equation}
can be written as
\begin{align}
  u_{2n+1} - u_{2n} &= \dt_{2n} D(u_{2n}) \\
  u_{2n+2} - u_{2n+1} &= - (u_{2n+1} - u_{2n}) + 2 \dt_{2n+1} D(u_{2n+1}),
\end{align}
where we have renumbered the steps so that the even numbered ones are the
results of complete RK steps and $\dt_n = \dt/2$.

For local time-stepping, the derivative values can depend on an additional set
of values $v_j$ (which have their own, similar, update equation), but where we
still expect the update rule to have the form of a linear combination:
\begin{equation}
  u_{n+1} - u_n = \sum_i A_n^i (u_n - u_{n-i})
  + \dt_n\, \sum_{i,j} B_n^{ij} D(u_{n-i}, v_{m(n)-j}).
\end{equation}

We now perform an element splitting as in Section~\ref{sec:element-splitting}
by writing $D(u, v) = V(u) + B(u, v)$.  Substituting this in gives
\begin{multline}\label{eq:general-framework:GL-split}
  u_{n+1} - u_n
  = \Big[\sum_i A_n^i (u_n - u_{n-i})
    + \dt_n\, \sum_i \Big(\sum_j B_n^{ij}\Big) V(u_{n-i})\Big] \\
  + \dt_n\, \sum_{i,j} B_n^{ij} B(u_{n-i}, v_{m(n)-j}).
\end{multline}
Since a general method must be independent of the details of the $V$ and $B$
functions, the bracketed terms in~\eqref{eq:general-framework:GL-split} must be
the standard GTS method operating with only the ``volume'' portion of the
equations, and the last term is a coupling correction.  Notably, the coupling
term does not require the function values directly, but only the value of the
coupling evaluated at those values.

When using a GTS method in an LTS framework, the~$u$ and $v$ will be evaluated
at the same sequence of times and the coefficients~$B_n^{ij}$ will be diagonal
in~$i,j$.  Comparing~\eqref{eq:general-framework:gts} to the bracketed term in~\eqref{eq:general-framework:GL-split}, we see that $\sum_j B_n^{ij} = B_n^i$,
so for a GTS method $B_n^{ij} = \delta_{ij} B_n^i$.  Combining all this, we
find that
\begin{equation}
  u_{n+1} - u_n
  = \Big[\sum_i A_n^i (u_n - u_{n-i})
    + \dt_n\, \sum_i B_n^i V(u_{n-i})\Big]
  + \dt_n\, \sum_i B_n^i B(u_{n-i}, v_{n-i}),
\end{equation}
that is, when using an arbitrary GTS integrator in a framework designed for
LTS, one can evaluate the volume term using the standard GTS rule and the
coupling contribution by using the usual update formula but with all the
non-derivative terms set to zero.  For the midpoint Runge-Kutta scheme above,
this gives the split rule
\begin{align}
  u_{2n+1} - u_{2n} &= \dt_{2n} V(u_{2n}) + \dt_{2n} B(u_{2n}, v_{2n}) \\
  u_{2n+2} - u_{2n+1} &=
  \Big[(u_{2n+1} - u_{2n+1-1}) + \dt_{2n+1} V(u_{2n+1})\Big]
  + \dt_{2n+1}\, B(u_{2n+1}, v_{2n+1}).
\end{align}

\section{Tables of coefficients for $2:1$ LTS rules}
\label{sec:tables}

Below are tables of coefficients
$a^s_{m;q^Aq^B}$ in~\eqref{eq:lts-form}
for order 2, 3, and 4 LTS rules with $2:1$
stepping, as well as the coefficients for transitioning between LTS and GTS
stepping in these cases.  The step patterns corresponding to these tables are
shown in Figure~\ref{fig:table-diagram}.
Additional values can be obtained using~\eqref{eq:lts-coef-small} and
\eqref{eq:lts-coef-large}.

\clearpage

For the transition rules, the number of steps requiring special coefficients
depends on the order of the integrator.  Only tables for steps affected by the
transition are shown below, after which either the $2:1$ rule or the GTS rule
should be used, as appropriate.

\begin{figure}
  \centering%
  \newcommand\hgap{5pt}
  \begin{tikzpicture}
    \Astep{-2}{-2\dt^A}
    \Astep{-1}{-\dt^A}
    \Astep{0}{0}
    \Astep{1}{\dt^A}
    \Astep{2}{2\dt^A}
    \Bstep{-4/2}{-4\dt^B}
    \Bstep{-3/2}{-3\dt^B}
    \Bstep{-2/2}{-2\dt^B}
    \Bstep{-1/2}{-\dt^B}
    \Bstep{0}{0}
    \Bstep{1/2}{\dt^B}
    \Bstep{2/2}{2\dt^B}
    \Bstep{3/2}{3\dt^B}
    \Bstep{4/2}{4\dt^B}
    \chartdots{5/2}
    \chartdots{-5/2}
    \Achart{1/2}{(a)}
    \Bchart{1/4}{(b)}
    \Bchart{3/4}{(c)}
  \end{tikzpicture}%
  \hspace{\hgap}%
  \begin{tikzpicture}
    \Astep{-2}{-2\dt^A}
    \Astep{-1}{-\dt^A}
    \Astep{0}{0}
    \Astep{1}{\dt^A}
    \Astep{2}{2\dt^A}
    \Bstep{-4/2}{-2\dt^A}
    \Bstep{-2/2}{-\dt^A}
    \Bstep{0}{0}
    \Bstep{1/2}{\dt^B}
    \Bstep{2/2}{2\dt^B}
    \Bstep{3/2}{3\dt^B}
    \Bstep{4/2}{4\dt^B}
    \chartdots{5/2}
    \chartdots{-5/2}
    \Achart{1/2}{(d0)}
    \Bchart{1/4}{(e0)}
    \Bchart{3/4}{(f0)}
    \Achart{3/2}{(d1)}
    \Bchart{5/4}{(e1)}
    \Bchart{7/4}{(f1)}
  \end{tikzpicture}%
  \hspace{\hgap}%
  \begin{tikzpicture}
    \Astep{-4/2}{-4\dt^B}
    \Astep{-3/2}{-3\dt^B}
    \Astep{-2/2}{-2\dt^B}
    \Astep{-1/2}{-\dt^B}
    \Astep{0}{0}
    \Astep{1}{\dt^A}
    \Astep{2}{2\dt^A}
    \Bstep{-4/2}{-4\dt^B}
    \Bstep{-3/2}{-3\dt^B}
    \Bstep{-2/2}{-2\dt^B}
    \Bstep{-1/2}{-\dt^B}
    \Bstep{0}{0}
    \Bstep{1/2}{\dt^B}
    \Bstep{2/2}{2\dt^B}
    \Bstep{3/2}{3\dt^B}
    \Bstep{4/2}{4\dt^B}
    \chartdots{5/2}
    \chartdots{-5/2}
    \Achart{1/2}{(g0)}
    \Bchart{1/4}{(h0)}
    \Bchart{3/4}{(i0)}
    \Achart{3/2}{(g1)}
    \Bchart{5/4}{(h1)}
    \Bchart{7/4}{(i1)}
  \end{tikzpicture}%

  \begin{tikzpicture}
    \Astep{-2}{-2\dt^A}
    \Astep{-1}{-\dt^A}
    \Astep{0}{0}
    \Astep{1/2}{\dt^B}
    \Astep{2/2}{2\dt^B}
    \Astep{3/2}{3\dt^B}
    \Astep{4/2}{4\dt^B}
    \Bstep{-4/2}{-4\dt^B}
    \Bstep{-3/2}{-3\dt^B}
    \Bstep{-2/2}{-2\dt^B}
    \Bstep{-1/2}{-\dt^B}
    \Bstep{0}{0}
    \Bstep{1/2}{\dt^B}
    \Bstep{2/2}{2\dt^B}
    \Bstep{3/2}{3\dt^B}
    \Bstep{4/2}{4\dt^B}
    \chartdots{5/2}
    \chartdots{-5/2}
    \Achart{1/4}{(j0)}
    \Bchart{1/4}{(j0)}
    \Achart{3/4}{(j1)}
    \Bchart{3/4}{(j1)}
    \Achart{5/4}{(j2)}
    \Bchart{5/4}{(j2)}
    \Achart{7/4}{(j3)}
    \Bchart{7/4}{(j3)}
  \end{tikzpicture}%
  \hspace{\hgap}%
  \begin{tikzpicture}
    \Astep{-2}{-2\dt^A}
    \Astep{-1}{-\dt^A}
    \Astep{0}{0}
    \Astep{1}{\dt^A}
    \Astep{2}{2\dt^A}
    \Bstep{-4/2}{-4\dt^B}
    \Bstep{-3/2}{-3\dt^B}
    \Bstep{-2/2}{-2\dt^B}
    \Bstep{-1/2}{-\dt^B}
    \Bstep{0}{0}
    \Bstep{1}{\dt^A}
    \Bstep{2}{2\dt^A}
    \chartdots{5/2}
    \chartdots{-5/2}
    \Achart{1/2}{(k0)}
    \Bchart{1/2}{(k0)}
    \Achart{3/2}{(k1)}
    \Bchart{3/2}{(k1)}
  \end{tikzpicture}%

  \caption{Step patterns during (a--c) steady state $2:1$ evolution, (d--f)
    transition to LTS by decreasing a step size, (g--i) transition to LTS by
    increasing a step size, (j) transition back to GTS by decreasing a step
    size, and (k) transition back to GTS by increasing a step size.  The labels
    correspond to the tables given in Appendix~\ref{sec:tables}.  The
    coefficients for transitioning back to GTS are the same for both elements.
    The numbered labels are extended upwards as necessary until the
    steady-state values are reached.}
  \label{fig:table-diagram}
\end{figure}

\input{lts_tables}

\bibliographystyle{siamplain}
\bibliography{references}
\end{document}

%% file: plots/step_pattern_burgers.tex
\begingroup
  \makeatletter
  \providecommand\color[2][]{%
    \GenericError{(gnuplot) \space\space\space\@spaces}{%
      Package color not loaded in conjunction with
      terminal option `colourtext'%
    }{See the gnuplot documentation for explanation.%
    }{Either use 'blacktext' in gnuplot or load the package
      color.sty in LaTeX.}%
    \renewcommand\color[2][]{}%
  }%
  \providecommand\includegraphics[2][]{%
    \GenericError{(gnuplot) \space\space\space\@spaces}{%
      Package graphicx or graphics not loaded%
    }{See the gnuplot documentation for explanation.%
    }{The gnuplot epslatex terminal needs graphicx.sty or graphics.sty.}%
    \renewcommand\includegraphics[2][]{}%
  }%
  \providecommand\rotatebox[2]{#2}%
  \@ifundefined{ifGPcolor}{%
    \newif\ifGPcolor
    \GPcolorfalse
  }{}%
  \@ifundefined{ifGPblacktext}{%
    \newif\ifGPblacktext
    \GPblacktexttrue
  }{}%
  \let\gplgaddtomacro\g@addto@macro
  \gdef\gplbacktext{}%
  \gdef\gplfronttext{}%
  \makeatother
  \ifGPblacktext
    \def\colorrgb#1{}%
    \def\colorgray#1{}%
  \else
    \ifGPcolor
      \def\colorrgb#1{\color[rgb]{#1}}%
      \def\colorgray#1{\color[gray]{#1}}%
      \expandafter\def\csname LTw\endcsname{\color{white}}%
      \expandafter\def\csname LTb\endcsname{\color{black}}%
      \expandafter\def\csname LTa\endcsname{\color{black}}%
      \expandafter\def\csname LT0\endcsname{\color[rgb]{1,0,0}}%
      \expandafter\def\csname LT1\endcsname{\color[rgb]{0,1,0}}%
      \expandafter\def\csname LT2\endcsname{\color[rgb]{0,0,1}}%
      \expandafter\def\csname LT3\endcsname{\color[rgb]{1,0,1}}%
      \expandafter\def\csname LT4\endcsname{\color[rgb]{0,1,1}}%
      \expandafter\def\csname LT5\endcsname{\color[rgb]{1,1,0}}%
      \expandafter\def\csname LT6\endcsname{\color[rgb]{0,0,0}}%
      \expandafter\def\csname LT7\endcsname{\color[rgb]{1,0.3,0}}%
      \expandafter\def\csname LT8\endcsname{\color[rgb]{0.5,0.5,0.5}}%
    \else
      \def\colorrgb#1{\color{black}}%
      \def\colorgray#1{\color[gray]{#1}}%
      \expandafter\def\csname LTw\endcsname{\color{white}}%
      \expandafter\def\csname LTb\endcsname{\color{black}}%
      \expandafter\def\csname LTa\endcsname{\color{black}}%
      \expandafter\def\csname LT0\endcsname{\color{black}}%
      \expandafter\def\csname LT1\endcsname{\color{black}}%
      \expandafter\def\csname LT2\endcsname{\color{black}}%
      \expandafter\def\csname LT3\endcsname{\color{black}}%
      \expandafter\def\csname LT4\endcsname{\color{black}}%
      \expandafter\def\csname LT5\endcsname{\color{black}}%
      \expandafter\def\csname LT6\endcsname{\color{black}}%
      \expandafter\def\csname LT7\endcsname{\color{black}}%
      \expandafter\def\csname LT8\endcsname{\color{black}}%
    \fi
  \fi
    \setlength{\unitlength}{0.0500bp}%
    \ifx\gptboxheight\undefined%
      \newlength{\gptboxheight}%
      \newlength{\gptboxwidth}%
      \newsavebox{\gptboxtext}%
    \fi%
    \setlength{\fboxrule}{0.5pt}%
    \setlength{\fboxsep}{1pt}%
\begin{picture}(7200.00,5040.00)%
    \gplgaddtomacro\gplbacktext{%
      \csname LTb\endcsname
      \put(814,1021){\makebox(0,0)[r]{\strut{}$0$}}%
      \put(814,1527){\makebox(0,0)[r]{\strut{}$0.2$}}%
      \put(814,2033){\makebox(0,0)[r]{\strut{}$0.4$}}%
      \put(814,2540){\makebox(0,0)[r]{\strut{}$0.6$}}%
      \put(814,3046){\makebox(0,0)[r]{\strut{}$0.8$}}%
      \put(814,3553){\makebox(0,0)[r]{\strut{}$1$}}%
      \put(814,4059){\makebox(0,0)[r]{\strut{}$1.2$}}%
      \put(814,4566){\makebox(0,0)[r]{\strut{}$1.4$}}%
      \put(1098,484){\makebox(0,0){\strut{}$0$}}%
      \put(1402,484){\makebox(0,0){\strut{}$1$}}%
      \put(1706,484){\makebox(0,0){\strut{}$2$}}%
      \put(2010,484){\makebox(0,0){\strut{}$3$}}%
      \put(2313,484){\makebox(0,0){\strut{}$4$}}%
      \put(2617,484){\makebox(0,0){\strut{}$5$}}%
      \put(2921,484){\makebox(0,0){\strut{}$6$}}%
      \put(3225,484){\makebox(0,0){\strut{}$7$}}%
      \put(3529,484){\makebox(0,0){\strut{}$8$}}%
      \put(3833,484){\makebox(0,0){\strut{}$9$}}%
      \put(4137,484){\makebox(0,0){\strut{}$10$}}%
      \put(4441,484){\makebox(0,0){\strut{}$11$}}%
      \put(4744,484){\makebox(0,0){\strut{}$12$}}%
      \put(5048,484){\makebox(0,0){\strut{}$13$}}%
      \put(5352,484){\makebox(0,0){\strut{}$14$}}%
      \put(5656,484){\makebox(0,0){\strut{}$15$}}%
    }%
    \gplgaddtomacro\gplfronttext{%
      \csname LTb\endcsname
      \put(308,2761){\rotatebox{-270}{\makebox(0,0){Time}}}%
      \put(3377,154){\makebox(0,0){Element}}%
      \csname LTb\endcsname
      \put(6305,1046){\makebox(0,0)[l]{$2^{-12}$}}%
      \put(6305,1732){\makebox(0,0)[l]{$2^{-11}$}}%
      \put(6305,2418){\makebox(0,0)[l]{$2^{-10}$}}%
      \put(6305,3104){\makebox(0,0)[l]{$2^{-9}$}}%
      \put(6305,3790){\makebox(0,0)[l]{$2^{-8}$}}%
      \put(6305,4476){\makebox(0,0)[l]{$2^{-7}$}}%
      \put(6899,2761){\rotatebox{-270}{\makebox(0,0){Step size}}}%
    }%
    \gplbacktext
    \put(0,0){\includegraphics{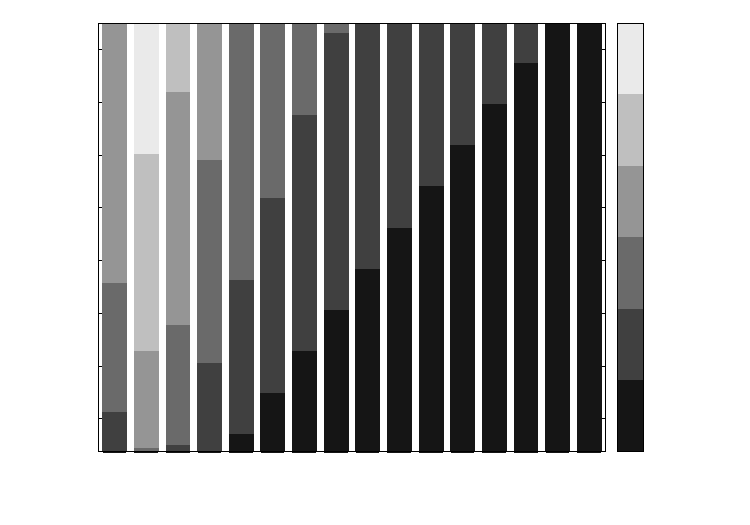}}%
    \gplfronttext
  \end{picture}%
\endgroup

%% file: plots/convergence.tex
\begingroup
  \makeatletter
  \providecommand\color[2][]{%
    \GenericError{(gnuplot) \space\space\space\@spaces}{%
      Package color not loaded in conjunction with
      terminal option `colourtext'%
    }{See the gnuplot documentation for explanation.%
    }{Either use 'blacktext' in gnuplot or load the package
      color.sty in LaTeX.}%
    \renewcommand\color[2][]{}%
  }%
  \providecommand\includegraphics[2][]{%
    \GenericError{(gnuplot) \space\space\space\@spaces}{%
      Package graphicx or graphics not loaded%
    }{See the gnuplot documentation for explanation.%
    }{The gnuplot epslatex terminal needs graphicx.sty or graphics.sty.}%
    \renewcommand\includegraphics[2][]{}%
  }%
  \providecommand\rotatebox[2]{#2}%
  \@ifundefined{ifGPcolor}{%
    \newif\ifGPcolor
    \GPcolorfalse
  }{}%
  \@ifundefined{ifGPblacktext}{%
    \newif\ifGPblacktext
    \GPblacktexttrue
  }{}%
  \let\gplgaddtomacro\g@addto@macro
  \gdef\gplbacktext{}%
  \gdef\gplfronttext{}%
  \makeatother
  \ifGPblacktext
    \def\colorrgb#1{}%
    \def\colorgray#1{}%
  \else
    \ifGPcolor
      \def\colorrgb#1{\color[rgb]{#1}}%
      \def\colorgray#1{\color[gray]{#1}}%
      \expandafter\def\csname LTw\endcsname{\color{white}}%
      \expandafter\def\csname LTb\endcsname{\color{black}}%
      \expandafter\def\csname LTa\endcsname{\color{black}}%
      \expandafter\def\csname LT0\endcsname{\color[rgb]{1,0,0}}%
      \expandafter\def\csname LT1\endcsname{\color[rgb]{0,1,0}}%
      \expandafter\def\csname LT2\endcsname{\color[rgb]{0,0,1}}%
      \expandafter\def\csname LT3\endcsname{\color[rgb]{1,0,1}}%
      \expandafter\def\csname LT4\endcsname{\color[rgb]{0,1,1}}%
      \expandafter\def\csname LT5\endcsname{\color[rgb]{1,1,0}}%
      \expandafter\def\csname LT6\endcsname{\color[rgb]{0,0,0}}%
      \expandafter\def\csname LT7\endcsname{\color[rgb]{1,0.3,0}}%
      \expandafter\def\csname LT8\endcsname{\color[rgb]{0.5,0.5,0.5}}%
    \else
      \def\colorrgb#1{\color{black}}%
      \def\colorgray#1{\color[gray]{#1}}%
      \expandafter\def\csname LTw\endcsname{\color{white}}%
      \expandafter\def\csname LTb\endcsname{\color{black}}%
      \expandafter\def\csname LTa\endcsname{\color{black}}%
      \expandafter\def\csname LT0\endcsname{\color{black}}%
      \expandafter\def\csname LT1\endcsname{\color{black}}%
      \expandafter\def\csname LT2\endcsname{\color{black}}%
      \expandafter\def\csname LT3\endcsname{\color{black}}%
      \expandafter\def\csname LT4\endcsname{\color{black}}%
      \expandafter\def\csname LT5\endcsname{\color{black}}%
      \expandafter\def\csname LT6\endcsname{\color{black}}%
      \expandafter\def\csname LT7\endcsname{\color{black}}%
      \expandafter\def\csname LT8\endcsname{\color{black}}%
    \fi
  \fi
    \setlength{\unitlength}{0.0500bp}%
    \ifx\gptboxheight\undefined%
      \newlength{\gptboxheight}%
      \newlength{\gptboxwidth}%
      \newsavebox{\gptboxtext}%
    \fi%
    \setlength{\fboxrule}{0.5pt}%
    \setlength{\fboxsep}{1pt}%
\begin{picture}(7200.00,5040.00)%
    \gplgaddtomacro\gplbacktext{%
      \csname LTb\endcsname
      \put(1078,704){\makebox(0,0)[r]{\strut{}$10^{-15}$}}%
      \put(1078,1161){\makebox(0,0)[r]{\strut{}$10^{-14}$}}%
      \put(1078,1618){\makebox(0,0)[r]{\strut{}$10^{-13}$}}%
      \put(1078,2076){\makebox(0,0)[r]{\strut{}$10^{-12}$}}%
      \put(1078,2533){\makebox(0,0)[r]{\strut{}$10^{-11}$}}%
      \put(1078,2990){\makebox(0,0)[r]{\strut{}$10^{-10}$}}%
      \put(1078,3447){\makebox(0,0)[r]{\strut{}$10^{-9}$}}%
      \put(1078,3905){\makebox(0,0)[r]{\strut{}$10^{-8}$}}%
      \put(1078,4362){\makebox(0,0)[r]{\strut{}$10^{-7}$}}%
      \put(1078,4819){\makebox(0,0)[r]{\strut{}$10^{-6}$}}%
      \put(1210,484){\makebox(0,0){\strut{}$10^{-5}$}}%
      \put(4006,484){\makebox(0,0){\strut{}$10^{-4}$}}%
      \put(6803,484){\makebox(0,0){\strut{}$10^{-3}$}}%
    }%
    \gplgaddtomacro\gplfronttext{%
      \csname LTb\endcsname
      \put(308,2761){\rotatebox{-270}{\makebox(0,0){$L_{\infty}$ error}}}%
      \put(4006,154){\makebox(0,0){Maximum allowed $u \dt$}}%
      \put(2197,4646){\makebox(0,0)[l]{\strut{}Order 4}}%
      \put(2197,4426){\makebox(0,0)[l]{\strut{}Order 5}}%
      \put(2197,4206){\makebox(0,0)[l]{\strut{}Order 6}}%
    }%
    \gplbacktext
    \put(0,0){\includegraphics{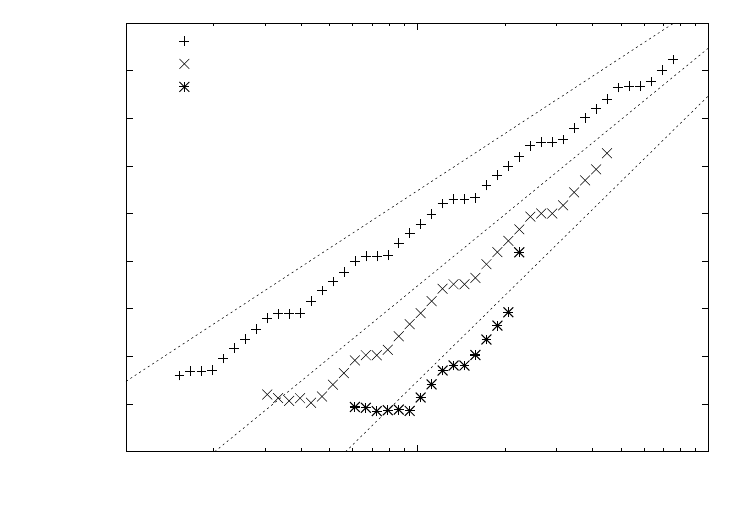}}%
    \gplfronttext
  \end{picture}%
\endgroup

%% file: plots/wave.tex
\begingroup
  \makeatletter
  \providecommand\color[2][]{%
    \GenericError{(gnuplot) \space\space\space\@spaces}{%
      Package color not loaded in conjunction with
      terminal option `colourtext'%
    }{See the gnuplot documentation for explanation.%
    }{Either use 'blacktext' in gnuplot or load the package
      color.sty in LaTeX.}%
    \renewcommand\color[2][]{}%
  }%
  \providecommand\includegraphics[2][]{%
    \GenericError{(gnuplot) \space\space\space\@spaces}{%
      Package graphicx or graphics not loaded%
    }{See the gnuplot documentation for explanation.%
    }{The gnuplot epslatex terminal needs graphicx.sty or graphics.sty.}%
    \renewcommand\includegraphics[2][]{}%
  }%
  \providecommand\rotatebox[2]{#2}%
  \@ifundefined{ifGPcolor}{%
    \newif\ifGPcolor
    \GPcolorfalse
  }{}%
  \@ifundefined{ifGPblacktext}{%
    \newif\ifGPblacktext
    \GPblacktexttrue
  }{}%
  \let\gplgaddtomacro\g@addto@macro
  \gdef\gplbacktext{}%
  \gdef\gplfronttext{}%
  \makeatother
  \ifGPblacktext
    \def\colorrgb#1{}%
    \def\colorgray#1{}%
  \else
    \ifGPcolor
      \def\colorrgb#1{\color[rgb]{#1}}%
      \def\colorgray#1{\color[gray]{#1}}%
      \expandafter\def\csname LTw\endcsname{\color{white}}%
      \expandafter\def\csname LTb\endcsname{\color{black}}%
      \expandafter\def\csname LTa\endcsname{\color{black}}%
      \expandafter\def\csname LT0\endcsname{\color[rgb]{1,0,0}}%
      \expandafter\def\csname LT1\endcsname{\color[rgb]{0,1,0}}%
      \expandafter\def\csname LT2\endcsname{\color[rgb]{0,0,1}}%
      \expandafter\def\csname LT3\endcsname{\color[rgb]{1,0,1}}%
      \expandafter\def\csname LT4\endcsname{\color[rgb]{0,1,1}}%
      \expandafter\def\csname LT5\endcsname{\color[rgb]{1,1,0}}%
      \expandafter\def\csname LT6\endcsname{\color[rgb]{0,0,0}}%
      \expandafter\def\csname LT7\endcsname{\color[rgb]{1,0.3,0}}%
      \expandafter\def\csname LT8\endcsname{\color[rgb]{0.5,0.5,0.5}}%
    \else
      \def\colorrgb#1{\color{black}}%
      \def\colorgray#1{\color[gray]{#1}}%
      \expandafter\def\csname LTw\endcsname{\color{white}}%
      \expandafter\def\csname LTb\endcsname{\color{black}}%
      \expandafter\def\csname LTa\endcsname{\color{black}}%
      \expandafter\def\csname LT0\endcsname{\color{black}}%
      \expandafter\def\csname LT1\endcsname{\color{black}}%
      \expandafter\def\csname LT2\endcsname{\color{black}}%
      \expandafter\def\csname LT3\endcsname{\color{black}}%
      \expandafter\def\csname LT4\endcsname{\color{black}}%
      \expandafter\def\csname LT5\endcsname{\color{black}}%
      \expandafter\def\csname LT6\endcsname{\color{black}}%
      \expandafter\def\csname LT7\endcsname{\color{black}}%
      \expandafter\def\csname LT8\endcsname{\color{black}}%
    \fi
  \fi
    \setlength{\unitlength}{0.0500bp}%
    \ifx\gptboxheight\undefined%
      \newlength{\gptboxheight}%
      \newlength{\gptboxwidth}%
      \newsavebox{\gptboxtext}%
    \fi%
    \setlength{\fboxrule}{0.5pt}%
    \setlength{\fboxsep}{1pt}%
\begin{picture}(7200.00,5040.00)%
    \gplgaddtomacro\gplbacktext{%
      \csname LTb\endcsname
      \put(946,704){\makebox(0,0)[r]{$-0.2$}}%
      \put(946,1218){\makebox(0,0)[r]{$0.0$}}%
      \put(946,1733){\makebox(0,0)[r]{$0.2$}}%
      \put(946,2247){\makebox(0,0)[r]{$0.4$}}%
      \put(946,2762){\makebox(0,0)[r]{$0.6$}}%
      \put(946,3276){\makebox(0,0)[r]{$0.8$}}%
      \put(946,3790){\makebox(0,0)[r]{$1.0$}}%
      \put(946,4305){\makebox(0,0)[r]{$1.2$}}%
      \put(946,4819){\makebox(0,0)[r]{$1.4$}}%
      \put(1650,484){\makebox(0,0){$-1.0$}}%
      \put(2566,484){\makebox(0,0){$-0.8$}}%
      \put(3482,484){\makebox(0,0){$-0.6$}}%
      \put(4398,484){\makebox(0,0){$-0.4$}}%
      \put(5314,484){\makebox(0,0){$-0.2$}}%
      \put(6230,484){\makebox(0,0){$0.0$}}%
    }%
    \gplgaddtomacro\gplfronttext{%
      \csname LTb\endcsname
      \put(308,2761){\rotatebox{-270}{\makebox(0,0){$u$}}}%
      \put(3940,154){\makebox(0,0){$x$}}%
      \put(5816,4646){\makebox(0,0)[r]{$t=0.0$}}%
      \put(5816,4426){\makebox(0,0)[r]{$t=0.1$}}%
      \put(5816,4206){\makebox(0,0)[r]{$t=0.2$}}%
      \put(5816,3986){\makebox(0,0)[r]{$t=0.3$}}%
      \put(5816,3766){\makebox(0,0)[r]{$t=0.4$}}%
      \put(5816,3546){\makebox(0,0)[r]{$t=0.5$}}%
    }%
    \gplbacktext
    \put(0,0){\includegraphics{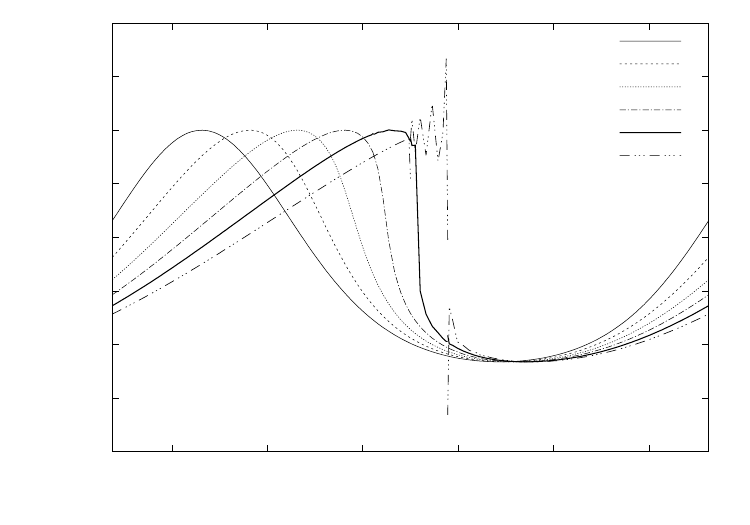}}%
    \gplfronttext
  \end{picture}%
\endgroup

%% file: plots/step_pattern_periodic.tex
\begingroup
  \makeatletter
  \providecommand\color[2][]{%
    \GenericError{(gnuplot) \space\space\space\@spaces}{%
      Package color not loaded in conjunction with
      terminal option `colourtext'%
    }{See the gnuplot documentation for explanation.%
    }{Either use 'blacktext' in gnuplot or load the package
      color.sty in LaTeX.}%
    \renewcommand\color[2][]{}%
  }%
  \providecommand\includegraphics[2][]{%
    \GenericError{(gnuplot) \space\space\space\@spaces}{%
      Package graphicx or graphics not loaded%
    }{See the gnuplot documentation for explanation.%
    }{The gnuplot epslatex terminal needs graphicx.sty or graphics.sty.}%
    \renewcommand\includegraphics[2][]{}%
  }%
  \providecommand\rotatebox[2]{#2}%
  \@ifundefined{ifGPcolor}{%
    \newif\ifGPcolor
    \GPcolorfalse
  }{}%
  \@ifundefined{ifGPblacktext}{%
    \newif\ifGPblacktext
    \GPblacktexttrue
  }{}%
  \let\gplgaddtomacro\g@addto@macro
  \gdef\gplbacktext{}%
  \gdef\gplfronttext{}%
  \makeatother
  \ifGPblacktext
    \def\colorrgb#1{}%
    \def\colorgray#1{}%
  \else
    \ifGPcolor
      \def\colorrgb#1{\color[rgb]{#1}}%
      \def\colorgray#1{\color[gray]{#1}}%
      \expandafter\def\csname LTw\endcsname{\color{white}}%
      \expandafter\def\csname LTb\endcsname{\color{black}}%
      \expandafter\def\csname LTa\endcsname{\color{black}}%
      \expandafter\def\csname LT0\endcsname{\color[rgb]{1,0,0}}%
      \expandafter\def\csname LT1\endcsname{\color[rgb]{0,1,0}}%
      \expandafter\def\csname LT2\endcsname{\color[rgb]{0,0,1}}%
      \expandafter\def\csname LT3\endcsname{\color[rgb]{1,0,1}}%
      \expandafter\def\csname LT4\endcsname{\color[rgb]{0,1,1}}%
      \expandafter\def\csname LT5\endcsname{\color[rgb]{1,1,0}}%
      \expandafter\def\csname LT6\endcsname{\color[rgb]{0,0,0}}%
      \expandafter\def\csname LT7\endcsname{\color[rgb]{1,0.3,0}}%
      \expandafter\def\csname LT8\endcsname{\color[rgb]{0.5,0.5,0.5}}%
    \else
      \def\colorrgb#1{\color{black}}%
      \def\colorgray#1{\color[gray]{#1}}%
      \expandafter\def\csname LTw\endcsname{\color{white}}%
      \expandafter\def\csname LTb\endcsname{\color{black}}%
      \expandafter\def\csname LTa\endcsname{\color{black}}%
      \expandafter\def\csname LT0\endcsname{\color{black}}%
      \expandafter\def\csname LT1\endcsname{\color{black}}%
      \expandafter\def\csname LT2\endcsname{\color{black}}%
      \expandafter\def\csname LT3\endcsname{\color{black}}%
      \expandafter\def\csname LT4\endcsname{\color{black}}%
      \expandafter\def\csname LT5\endcsname{\color{black}}%
      \expandafter\def\csname LT6\endcsname{\color{black}}%
      \expandafter\def\csname LT7\endcsname{\color{black}}%
      \expandafter\def\csname LT8\endcsname{\color{black}}%
    \fi
  \fi
    \setlength{\unitlength}{0.0500bp}%
    \ifx\gptboxheight\undefined%
      \newlength{\gptboxheight}%
      \newlength{\gptboxwidth}%
      \newsavebox{\gptboxtext}%
    \fi%
    \setlength{\fboxrule}{0.5pt}%
    \setlength{\fboxsep}{1pt}%
\begin{picture}(7200.00,5040.00)%
    \gplgaddtomacro\gplbacktext{%
      \csname LTb\endcsname
      \put(814,704){\makebox(0,0)[r]{\strut{}$0$}}%
      \put(814,1527){\makebox(0,0)[r]{\strut{}$0.1$}}%
      \put(814,2350){\makebox(0,0)[r]{\strut{}$0.2$}}%
      \put(814,3173){\makebox(0,0)[r]{\strut{}$0.3$}}%
      \put(814,3996){\makebox(0,0)[r]{\strut{}$0.4$}}%
      \put(814,4819){\makebox(0,0)[r]{\strut{}$0.5$}}%
      \put(1098,484){\makebox(0,0){\strut{}$0$}}%
      \put(1402,484){\makebox(0,0){\strut{}$1$}}%
      \put(1706,484){\makebox(0,0){\strut{}$2$}}%
      \put(2010,484){\makebox(0,0){\strut{}$3$}}%
      \put(2313,484){\makebox(0,0){\strut{}$4$}}%
      \put(2617,484){\makebox(0,0){\strut{}$5$}}%
      \put(2921,484){\makebox(0,0){\strut{}$6$}}%
      \put(3225,484){\makebox(0,0){\strut{}$7$}}%
      \put(3529,484){\makebox(0,0){\strut{}$8$}}%
      \put(3833,484){\makebox(0,0){\strut{}$9$}}%
      \put(4137,484){\makebox(0,0){\strut{}$10$}}%
      \put(4441,484){\makebox(0,0){\strut{}$11$}}%
      \put(4744,484){\makebox(0,0){\strut{}$12$}}%
      \put(5048,484){\makebox(0,0){\strut{}$13$}}%
      \put(5352,484){\makebox(0,0){\strut{}$14$}}%
      \put(5656,484){\makebox(0,0){\strut{}$15$}}%
    }%
    \gplgaddtomacro\gplfronttext{%
      \csname LTb\endcsname
      \put(308,2761){\rotatebox{-270}{\makebox(0,0){Time}}}%
      \put(3377,154){\makebox(0,0){Element}}%
      \csname LTb\endcsname
      \put(6305,1218){\makebox(0,0)[l]{$2^{-13}$}}%
      \put(6305,2247){\makebox(0,0)[l]{$2^{-12}$}}%
      \put(6305,3275){\makebox(0,0)[l]{$2^{-11}$}}%
      \put(6305,4304){\makebox(0,0)[l]{$2^{-10}$}}%
      \put(6899,2761){\rotatebox{-270}{\makebox(0,0){Step size}}}%
    }%
    \gplbacktext
    \put(0,0){\includegraphics{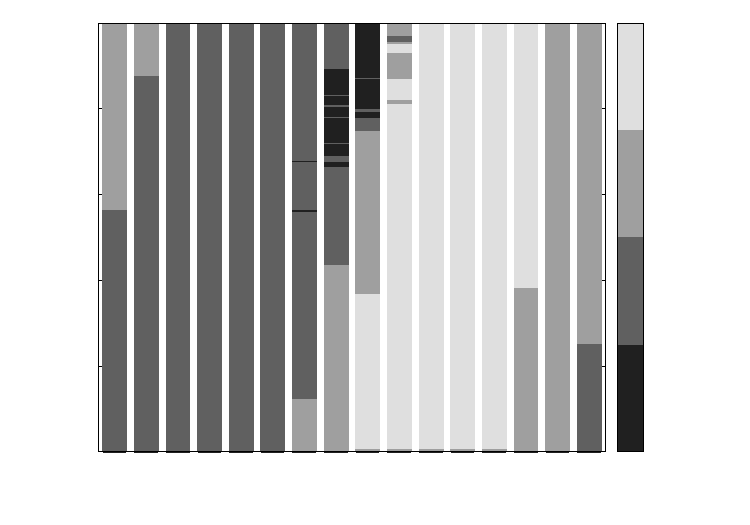}}%
    \gplfronttext
  \end{picture}%
\endgroup

%% file: plots/conservation.tex
\begingroup
  \makeatletter
  \providecommand\color[2][]{%
    \GenericError{(gnuplot) \space\space\space\@spaces}{%
      Package color not loaded in conjunction with
      terminal option `colourtext'%
    }{See the gnuplot documentation for explanation.%
    }{Either use 'blacktext' in gnuplot or load the package
      color.sty in LaTeX.}%
    \renewcommand\color[2][]{}%
  }%
  \providecommand\includegraphics[2][]{%
    \GenericError{(gnuplot) \space\space\space\@spaces}{%
      Package graphicx or graphics not loaded%
    }{See the gnuplot documentation for explanation.%
    }{The gnuplot epslatex terminal needs graphicx.sty or graphics.sty.}%
    \renewcommand\includegraphics[2][]{}%
  }%
  \providecommand\rotatebox[2]{#2}%
  \@ifundefined{ifGPcolor}{%
    \newif\ifGPcolor
    \GPcolorfalse
  }{}%
  \@ifundefined{ifGPblacktext}{%
    \newif\ifGPblacktext
    \GPblacktexttrue
  }{}%
  \let\gplgaddtomacro\g@addto@macro
  \gdef\gplbacktext{}%
  \gdef\gplfronttext{}%
  \makeatother
  \ifGPblacktext
    \def\colorrgb#1{}%
    \def\colorgray#1{}%
  \else
    \ifGPcolor
      \def\colorrgb#1{\color[rgb]{#1}}%
      \def\colorgray#1{\color[gray]{#1}}%
      \expandafter\def\csname LTw\endcsname{\color{white}}%
      \expandafter\def\csname LTb\endcsname{\color{black}}%
      \expandafter\def\csname LTa\endcsname{\color{black}}%
      \expandafter\def\csname LT0\endcsname{\color[rgb]{1,0,0}}%
      \expandafter\def\csname LT1\endcsname{\color[rgb]{0,1,0}}%
      \expandafter\def\csname LT2\endcsname{\color[rgb]{0,0,1}}%
      \expandafter\def\csname LT3\endcsname{\color[rgb]{1,0,1}}%
      \expandafter\def\csname LT4\endcsname{\color[rgb]{0,1,1}}%
      \expandafter\def\csname LT5\endcsname{\color[rgb]{1,1,0}}%
      \expandafter\def\csname LT6\endcsname{\color[rgb]{0,0,0}}%
      \expandafter\def\csname LT7\endcsname{\color[rgb]{1,0.3,0}}%
      \expandafter\def\csname LT8\endcsname{\color[rgb]{0.5,0.5,0.5}}%
    \else
      \def\colorrgb#1{\color{black}}%
      \def\colorgray#1{\color[gray]{#1}}%
      \expandafter\def\csname LTw\endcsname{\color{white}}%
      \expandafter\def\csname LTb\endcsname{\color{black}}%
      \expandafter\def\csname LTa\endcsname{\color{black}}%
      \expandafter\def\csname LT0\endcsname{\color{black}}%
      \expandafter\def\csname LT1\endcsname{\color{black}}%
      \expandafter\def\csname LT2\endcsname{\color{black}}%
      \expandafter\def\csname LT3\endcsname{\color{black}}%
      \expandafter\def\csname LT4\endcsname{\color{black}}%
      \expandafter\def\csname LT5\endcsname{\color{black}}%
      \expandafter\def\csname LT6\endcsname{\color{black}}%
      \expandafter\def\csname LT7\endcsname{\color{black}}%
      \expandafter\def\csname LT8\endcsname{\color{black}}%
    \fi
  \fi
    \setlength{\unitlength}{0.0500bp}%
    \ifx\gptboxheight\undefined%
      \newlength{\gptboxheight}%
      \newlength{\gptboxwidth}%
      \newsavebox{\gptboxtext}%
    \fi%
    \setlength{\fboxrule}{0.5pt}%
    \setlength{\fboxsep}{1pt}%
\begin{picture}(7200.00,5040.00)%
    \gplgaddtomacro\gplbacktext{%
      \csname LTb\endcsname
      \put(946,704){\makebox(0,0)[r]{$-2.5$}}%
      \put(946,1116){\makebox(0,0)[r]{$-2.0$}}%
      \put(946,1527){\makebox(0,0)[r]{$-1.5$}}%
      \put(946,1939){\makebox(0,0)[r]{$-1.0$}}%
      \put(946,2350){\makebox(0,0)[r]{$-0.5$}}%
      \put(946,2762){\makebox(0,0)[r]{$0.0$}}%
      \put(946,3173){\makebox(0,0)[r]{$0.5$}}%
      \put(946,3585){\makebox(0,0)[r]{$1.0$}}%
      \put(946,3996){\makebox(0,0)[r]{$1.5$}}%
      \put(946,4408){\makebox(0,0)[r]{$2.0$}}%
      \put(946,4819){\makebox(0,0)[r]{$2.5$}}%
      \put(1078,484){\makebox(0,0){\strut{}$0$}}%
      \put(2223,484){\makebox(0,0){\strut{}$200$}}%
      \put(3368,484){\makebox(0,0){\strut{}$400$}}%
      \put(4513,484){\makebox(0,0){\strut{}$600$}}%
      \put(5658,484){\makebox(0,0){\strut{}$800$}}%
      \put(6803,484){\makebox(0,0){\strut{}$1000$}}%
    }%
    \gplgaddtomacro\gplfronttext{%
      \csname LTb\endcsname
      \put(308,2761){\rotatebox{-270}{\makebox(0,0){$(\Delta \int\! u) / 10^{-14}$}}}%
      \put(3940,154){\makebox(0,0){Time}}%
    }%
    \gplbacktext
    \put(0,0){\includegraphics{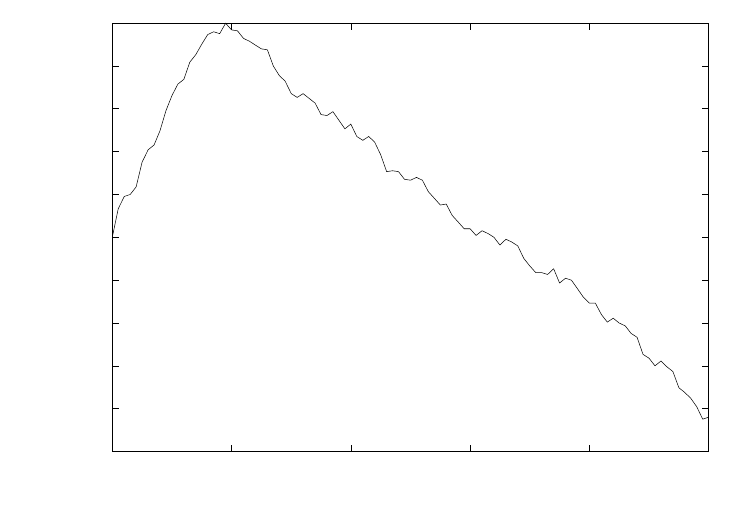}}%
    \gplfronttext
  \end{picture}%
\endgroup

%% file: plots/lts_cfl.tex
\begingroup
  \makeatletter
  \providecommand\color[2][]{%
    \GenericError{(gnuplot) \space\space\space\@spaces}{%
      Package color not loaded in conjunction with
      terminal option `colourtext'%
    }{See the gnuplot documentation for explanation.%
    }{Either use 'blacktext' in gnuplot or load the package
      color.sty in LaTeX.}%
    \renewcommand\color[2][]{}%
  }%
  \providecommand\includegraphics[2][]{%
    \GenericError{(gnuplot) \space\space\space\@spaces}{%
      Package graphicx or graphics not loaded%
    }{See the gnuplot documentation for explanation.%
    }{The gnuplot epslatex terminal needs graphicx.sty or graphics.sty.}%
    \renewcommand\includegraphics[2][]{}%
  }%
  \providecommand\rotatebox[2]{#2}%
  \@ifundefined{ifGPcolor}{%
    \newif\ifGPcolor
    \GPcolorfalse
  }{}%
  \@ifundefined{ifGPblacktext}{%
    \newif\ifGPblacktext
    \GPblacktexttrue
  }{}%
  \let\gplgaddtomacro\g@addto@macro
  \gdef\gplbacktext{}%
  \gdef\gplfronttext{}%
  \makeatother
  \ifGPblacktext
    \def\colorrgb#1{}%
    \def\colorgray#1{}%
  \else
    \ifGPcolor
      \def\colorrgb#1{\color[rgb]{#1}}%
      \def\colorgray#1{\color[gray]{#1}}%
      \expandafter\def\csname LTw\endcsname{\color{white}}%
      \expandafter\def\csname LTb\endcsname{\color{black}}%
      \expandafter\def\csname LTa\endcsname{\color{black}}%
      \expandafter\def\csname LT0\endcsname{\color[rgb]{1,0,0}}%
      \expandafter\def\csname LT1\endcsname{\color[rgb]{0,1,0}}%
      \expandafter\def\csname LT2\endcsname{\color[rgb]{0,0,1}}%
      \expandafter\def\csname LT3\endcsname{\color[rgb]{1,0,1}}%
      \expandafter\def\csname LT4\endcsname{\color[rgb]{0,1,1}}%
      \expandafter\def\csname LT5\endcsname{\color[rgb]{1,1,0}}%
      \expandafter\def\csname LT6\endcsname{\color[rgb]{0,0,0}}%
      \expandafter\def\csname LT7\endcsname{\color[rgb]{1,0.3,0}}%
      \expandafter\def\csname LT8\endcsname{\color[rgb]{0.5,0.5,0.5}}%
    \else
      \def\colorrgb#1{\color{black}}%
      \def\colorgray#1{\color[gray]{#1}}%
      \expandafter\def\csname LTw\endcsname{\color{white}}%
      \expandafter\def\csname LTb\endcsname{\color{black}}%
      \expandafter\def\csname LTa\endcsname{\color{black}}%
      \expandafter\def\csname LT0\endcsname{\color{black}}%
      \expandafter\def\csname LT1\endcsname{\color{black}}%
      \expandafter\def\csname LT2\endcsname{\color{black}}%
      \expandafter\def\csname LT3\endcsname{\color{black}}%
      \expandafter\def\csname LT4\endcsname{\color{black}}%
      \expandafter\def\csname LT5\endcsname{\color{black}}%
      \expandafter\def\csname LT6\endcsname{\color{black}}%
      \expandafter\def\csname LT7\endcsname{\color{black}}%
      \expandafter\def\csname LT8\endcsname{\color{black}}%
    \fi
  \fi
    \setlength{\unitlength}{0.0500bp}%
    \ifx\gptboxheight\undefined%
      \newlength{\gptboxheight}%
      \newlength{\gptboxwidth}%
      \newsavebox{\gptboxtext}%
    \fi%
    \setlength{\fboxrule}{0.5pt}%
    \setlength{\fboxsep}{1pt}%
\begin{picture}(7200.00,5040.00)%
    \gplgaddtomacro\gplbacktext{%
      \csname LTb\endcsname
      \put(946,704){\makebox(0,0)[r]{\strut{}0.70}}%
      \put(946,1390){\makebox(0,0)[r]{\strut{}0.75}}%
      \put(946,2076){\makebox(0,0)[r]{\strut{}0.80}}%
      \put(946,2762){\makebox(0,0)[r]{\strut{}0.85}}%
      \put(946,3447){\makebox(0,0)[r]{\strut{}0.90}}%
      \put(946,4133){\makebox(0,0)[r]{\strut{}0.95}}%
      \put(946,4819){\makebox(0,0)[r]{\strut{}1.00}}%
      \put(1078,484){\makebox(0,0){\strut{}$2$}}%
      \put(1794,484){\makebox(0,0){\strut{}$3$}}%
      \put(2509,484){\makebox(0,0){\strut{}$4$}}%
      \put(3225,484){\makebox(0,0){\strut{}$5$}}%
      \put(3941,484){\makebox(0,0){\strut{}$6$}}%
      \put(4656,484){\makebox(0,0){\strut{}$7$}}%
      \put(5372,484){\makebox(0,0){\strut{}$8$}}%
      \put(6087,484){\makebox(0,0){\strut{}$9$}}%
      \put(6803,484){\makebox(0,0){\strut{}$10$}}%
    }%
    \gplgaddtomacro\gplfronttext{%
      \csname LTb\endcsname
      \put(308,2761){\rotatebox{-270}{\makebox(0,0){$\displaystyle s = \frac{C}{C_{\text{advection}}}$}}}%
      \put(3940,154){\makebox(0,0){Higher refinement level $n$}}%
      \put(5816,4646){\makebox(0,0)[r]{\strut{}LTS}}%
      \put(5816,4426){\makebox(0,0)[r]{\strut{}GTS}}%
    }%
    \gplbacktext
    \put(0,0){\includegraphics{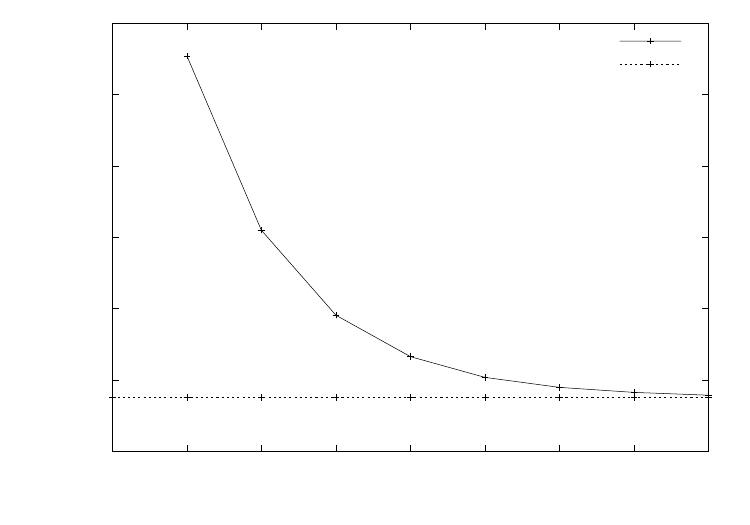}}%
    \gplfronttext
  \end{picture}%
\endgroup

%% file: plots/resolution_convergence.tex
\begingroup
  \makeatletter
  \providecommand\color[2][]{%
    \GenericError{(gnuplot) \space\space\space\@spaces}{%
      Package color not loaded in conjunction with
      terminal option `colourtext'%
    }{See the gnuplot documentation for explanation.%
    }{Either use 'blacktext' in gnuplot or load the package
      color.sty in LaTeX.}%
    \renewcommand\color[2][]{}%
  }%
  \providecommand\includegraphics[2][]{%
    \GenericError{(gnuplot) \space\space\space\@spaces}{%
      Package graphicx or graphics not loaded%
    }{See the gnuplot documentation for explanation.%
    }{The gnuplot epslatex terminal needs graphicx.sty or graphics.sty.}%
    \renewcommand\includegraphics[2][]{}%
  }%
  \providecommand\rotatebox[2]{#2}%
  \@ifundefined{ifGPcolor}{%
    \newif\ifGPcolor
    \GPcolorfalse
  }{}%
  \@ifundefined{ifGPblacktext}{%
    \newif\ifGPblacktext
    \GPblacktexttrue
  }{}%
  \let\gplgaddtomacro\g@addto@macro
  \gdef\gplbacktext{}%
  \gdef\gplfronttext{}%
  \makeatother
  \ifGPblacktext
    \def\colorrgb#1{}%
    \def\colorgray#1{}%
  \else
    \ifGPcolor
      \def\colorrgb#1{\color[rgb]{#1}}%
      \def\colorgray#1{\color[gray]{#1}}%
      \expandafter\def\csname LTw\endcsname{\color{white}}%
      \expandafter\def\csname LTb\endcsname{\color{black}}%
      \expandafter\def\csname LTa\endcsname{\color{black}}%
      \expandafter\def\csname LT0\endcsname{\color[rgb]{1,0,0}}%
      \expandafter\def\csname LT1\endcsname{\color[rgb]{0,1,0}}%
      \expandafter\def\csname LT2\endcsname{\color[rgb]{0,0,1}}%
      \expandafter\def\csname LT3\endcsname{\color[rgb]{1,0,1}}%
      \expandafter\def\csname LT4\endcsname{\color[rgb]{0,1,1}}%
      \expandafter\def\csname LT5\endcsname{\color[rgb]{1,1,0}}%
      \expandafter\def\csname LT6\endcsname{\color[rgb]{0,0,0}}%
      \expandafter\def\csname LT7\endcsname{\color[rgb]{1,0.3,0}}%
      \expandafter\def\csname LT8\endcsname{\color[rgb]{0.5,0.5,0.5}}%
    \else
      \def\colorrgb#1{\color{black}}%
      \def\colorgray#1{\color[gray]{#1}}%
      \expandafter\def\csname LTw\endcsname{\color{white}}%
      \expandafter\def\csname LTb\endcsname{\color{black}}%
      \expandafter\def\csname LTa\endcsname{\color{black}}%
      \expandafter\def\csname LT0\endcsname{\color{black}}%
      \expandafter\def\csname LT1\endcsname{\color{black}}%
      \expandafter\def\csname LT2\endcsname{\color{black}}%
      \expandafter\def\csname LT3\endcsname{\color{black}}%
      \expandafter\def\csname LT4\endcsname{\color{black}}%
      \expandafter\def\csname LT5\endcsname{\color{black}}%
      \expandafter\def\csname LT6\endcsname{\color{black}}%
      \expandafter\def\csname LT7\endcsname{\color{black}}%
      \expandafter\def\csname LT8\endcsname{\color{black}}%
    \fi
  \fi
    \setlength{\unitlength}{0.0500bp}%
    \ifx\gptboxheight\undefined%
      \newlength{\gptboxheight}%
      \newlength{\gptboxwidth}%
      \newsavebox{\gptboxtext}%
    \fi%
    \setlength{\fboxrule}{0.5pt}%
    \setlength{\fboxsep}{1pt}%
\begin{picture}(7200.00,5040.00)%
    \gplgaddtomacro\gplbacktext{%
      \csname LTb\endcsname
      \put(1078,704){\makebox(0,0)[r]{$10^{-14}$}}%
      \put(1078,1292){\makebox(0,0)[r]{$10^{-12}$}}%
      \put(1078,1880){\makebox(0,0)[r]{$10^{-10}$}}%
      \put(1078,2468){\makebox(0,0)[r]{$10^{-8}$}}%
      \put(1078,3055){\makebox(0,0)[r]{$10^{-6}$}}%
      \put(1078,3643){\makebox(0,0)[r]{$10^{-4}$}}%
      \put(1078,4231){\makebox(0,0)[r]{$10^{-2}$}}%
      \put(1078,4819){\makebox(0,0)[r]{$10^{0}$}}%
      \put(2238,484){\makebox(0,0){$10^{-2}$}}%
      \put(3786,484){\makebox(0,0){$10^{-1}$}}%
      \put(5334,484){\makebox(0,0){$10^{0}$}}%
      \put(5026,3028){\makebox(0,0){$\propto h^5$}}%
      \put(2374,1164){\makebox(0,0){$\propto h^3$}}%
      \put(2374,2225){\makebox(0,0){$\propto h^2$}}%
    }%
    \gplgaddtomacro\gplfronttext{%
      \csname LTb\endcsname
      \put(308,2761){\rotatebox{-270}{\makebox(0,0){$L_{\infty}$ error}}}%
      \put(4006,154){\makebox(0,0){Size of larger elements $h$}}%
      \put(2197,4646){\makebox(0,0)[l]{$A = 0.40$}}%
      \put(2197,4426){\makebox(0,0)[l]{$A = 0.41$}}%
    }%
    \gplbacktext
    \put(0,0){\includegraphics{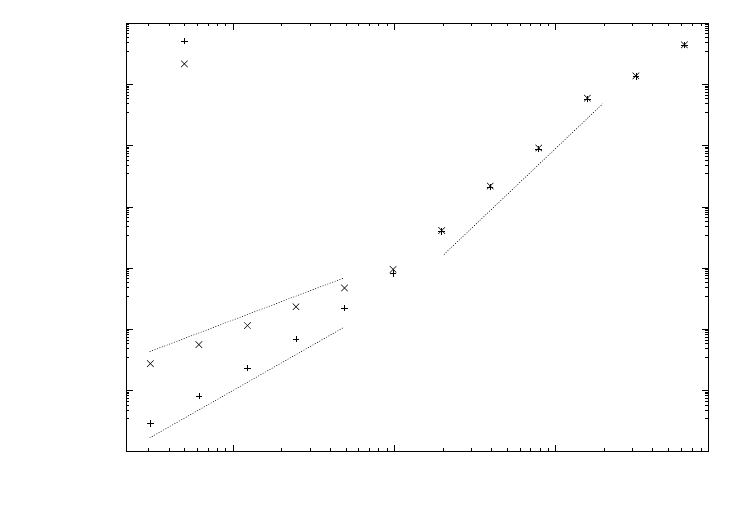}}%
    \gplfronttext
  \end{picture}%
\endgroup

%% file: lts_tables.tex
\begingroup

\newcommand\hgap{10pt}
\newcommand\vgap{10pt}

\subsection{Order 2}

\hspace{0pt} 

\begin{longtable}{s|sss@{}p{\hgap}@{}s|ss@{}p{\hgap}@{}s|ss}
\multicolumn{12}{c}{\textbf{LTS $2:1$ rule}} \\*
\multicolumn{1}{c|}{\text{(a)}} & \dt^B & 0 & -\dt^B && 
\multicolumn{1}{c|}{\text{(b)}} & 0 & -\dt^B && 
\multicolumn{1}{c|}{\text{(c)}} & \dt^B & 0 \\* 
\hhline{----~---~---}\nopagebreak
0 & \frac{9}{8} & \frac{1}{2} & -\frac{1}{8} && 
0 & \frac{3}{2} & -\frac{1}{4} && 
0 & \frac{9}{4} & -\frac{1}{2} \\* 
-\dt^A & -\frac{3}{8} & 0 & -\frac{1}{8} && 
-\dt^A & 0 & -\frac{1}{4} && 
-\dt^A & -\frac{3}{4} & 0 
\vspace{\vgap} \\
\multicolumn{12}{c}{\textbf{Transition to LTS by decreasing a step size}} \\*
\multicolumn{1}{c|}{\text{(d0)}} & \dt^B & 0 & -\dt^A && 
\multicolumn{1}{c|}{\text{(e0)}} & 0 & -\dt^A && 
\multicolumn{1}{c|}{\text{(f0)}} & \dt^B & 0 \\* 
\hhline{----~---~---}\nopagebreak
0 & \frac{9}{8} & \frac{3}{8} & 0 && 
0 & \frac{5}{4} & 0 && 
0 & \frac{9}{4} & -\frac{1}{2} \\* 
-\dt^A & -\frac{3}{8} & 0 & -\frac{1}{8} && 
-\dt^A & 0 & -\frac{1}{4} && 
-\dt^A & -\frac{3}{4} & 0 
\vspace{\vgap} \\
\multicolumn{12}{c}{\textbf{Transition to LTS by increasing a step size}} \\*
\multicolumn{1}{c|}{\text{(g0)}} & \dt^B & 0 & -\dt^B && 
\multicolumn{1}{c|}{\text{(h0)}} & 0 & -\dt^B && 
\multicolumn{1}{c|}{\text{(i0)}} & \dt^B & 0 \\* 
\hhline{----~---~---}\nopagebreak
0 & \frac{3}{2} & \frac{1}{2} & 0 && 
0 & \frac{3}{2} & 0 && 
0 & 3 & -\frac{1}{2} \\* 
-\dt^B & -\frac{3}{4} & 0 & -\frac{1}{4} && 
-\dt^B & 0 & -\frac{1}{2} && 
-\dt^B & -\frac{3}{2} & 0 
\vspace{\vgap} \\
\multicolumn{12}{c}{\textbf{Transitioning to GTS}} \\*
\multicolumn{1}{c|}{\text{(j0)}} & 0 & -\dt^B &&&\multicolumn{3}{c}{}&& 
\multicolumn{1}{c|}{\text{(k0)}} & 0 & -\dt^B \\* 
\hhline{---~~~~~~---}\nopagebreak
0 & \frac{3}{2} & -\frac{1}{4} &&&\multicolumn{3}{c}{}&& 
0 & 2 & -\frac{1}{2} \\* 
-\dt^A & 0 & -\frac{1}{4} &&&\multicolumn{3}{c}{}&& 
-\dt^A & 0 & -\frac{1}{2} 
\end{longtable}

\subsection{Order 3}

\hspace{0pt} 

\begin{longtable}{s|ssss@{}p{\hgap}@{}s|ssss}
\multicolumn{11}{c}{\textbf{LTS $2:1$ rule}} \\*
\multicolumn{1}{c|}{\text{(a)}} & \dt^B & 0 & -\dt^B & -2\dt^B && 
\multicolumn{1}{c|}{\text{(b)}} & 0 & -\dt^B & -2\dt^B \\* 
\hhline{-----~----}\nopagebreak
0 & \frac{115}{64} & \frac{7}{24} & -\frac{11}{64} & 0 && 
0 & \frac{23}{12} & -\frac{1}{2} & 0 \\* 
-\dt^A & -\frac{115}{96} & 0 & -\frac{11}{32} & \frac{5}{24} && 
-\dt^A & 0 & -1 & \frac{5}{12} \\* 
-2\dt^A & \frac{23}{64} & 0 & \frac{11}{192} & 0 && 
-2\dt^A & 0 & \frac{1}{6} & 0 
\vspace{\vgap} \\
\multicolumn{1}{c|}{\text{(c)}} & \dt^B & 0 & -\dt^B \\* 
\hhline{----}\nopagebreak
0 & \frac{115}{32} & -\frac{4}{3} & \frac{5}{32} \\* 
-\dt^A & -\frac{115}{48} & 0 & \frac{5}{16} \\* 
-2\dt^A & \frac{23}{32} & 0 & -\frac{5}{96} 
\vspace{\vgap} \\
\multicolumn{11}{c}{\textbf{Transition to LTS by decreasing a step size}} \\*
\multicolumn{1}{c|}{\text{(d0)}} & \dt^B & 0 & -\dt^A & -2\dt^A && 
\multicolumn{1}{c|}{\text{(e0)}} & 0 & -\dt^A & -2\dt^A \\* 
\hhline{-----~----}\nopagebreak
0 & \frac{5}{3} & \frac{1}{4} & 0 & 0 && 
0 & \frac{17}{12} & 0 & 0 \\* 
-\dt^A & -\frac{10}{9} & 0 & -\frac{2}{9} & 0 && 
-\dt^A & 0 & -\frac{7}{12} & 0 \\* 
-2\dt^A & \frac{1}{3} & 0 & 0 & \frac{1}{12} && 
-2\dt^A & 0 & 0 & \frac{1}{6} 
\vspace{\vgap} \\
\multicolumn{1}{c|}{\text{(f0)}} & \dt^B & 0 & -\dt^A \\* 
\hhline{----}\nopagebreak
0 & \frac{10}{3} & -\frac{11}{12} & 0 \\* 
-\dt^A & -\frac{20}{9} & 0 & \frac{5}{36} \\* 
-2\dt^A & \frac{2}{3} & 0 & 0 
\vspace{\vgap} \\
\multicolumn{11}{c}{\textbf{Transition to LTS by increasing a step size}} \\*
\multicolumn{1}{c|}{\text{(g0)}} & \dt^B & 0 & -\dt^B & -2\dt^B && 
\multicolumn{1}{c|}{\text{(h0)}} & 0 & -\dt^B & -2\dt^B \\* 
\hhline{-----~----}\nopagebreak
0 & \frac{23}{8} & \frac{7}{24} & 0 & 0 && 
0 & \frac{23}{12} & 0 & 0 \\* 
-\dt^B & -\frac{23}{8} & 0 & -\frac{11}{24} & 0 && 
-\dt^B & 0 & -\frac{4}{3} & 0 \\* 
-2\dt^B & \frac{23}{24} & 0 & 0 & \frac{5}{24} && 
-2\dt^B & 0 & 0 & \frac{5}{12} 
\vspace{\vgap} \\
\multicolumn{1}{c|}{\text{(i0)}} & \dt^B & 0 & -\dt^B &&& 
\multicolumn{1}{c|}{\text{(g1)}} & 3\dt^B & 2\dt^B & \dt^B & 0 \\* 
\hhline{----~~-----}\nopagebreak
0 & \frac{23}{4} & -\frac{4}{3} & 0 &&& 
\dt^A & \frac{23}{12} & \frac{7}{24} & -\frac{11}{72} & 0 \\* 
-\dt^B & -\frac{23}{4} & 0 & \frac{5}{12} &&& 
0 & -\frac{23}{12} & 0 & -\frac{11}{24} & \frac{5}{24} \\* 
-2\dt^B & \frac{23}{12} & 0 & 0 &&& 
-\dt^B & \frac{23}{24} & 0 & \frac{11}{72} & 0 
\vspace{\vgap} \\
\multicolumn{1}{c|}{\text{(h1)}} & 2\dt^B & \dt^B & 0 &&& 
\multicolumn{1}{c|}{\text{(i1)}} & 3\dt^B & 2\dt^B & \dt^B \\* 
\hhline{----~~----}\nopagebreak
\dt^A & \frac{23}{12} & -\frac{4}{9} & 0 &&& 
\dt^A & \frac{23}{6} & -\frac{4}{3} & \frac{5}{36} \\* 
0 & 0 & -\frac{4}{3} & \frac{5}{12} &&& 
0 & -\frac{23}{6} & 0 & \frac{5}{12} \\* 
-\dt^B & 0 & \frac{4}{9} & 0 &&& 
-\dt^B & \frac{23}{12} & 0 & -\frac{5}{36} 
\vspace{\vgap} \\
\multicolumn{11}{c}{\textbf{Transitioning to GTS by decreasing a step size}} \\*
\multicolumn{1}{c|}{\text{(j0)}} & 0 & -\dt^B & -2\dt^B &&& 
\multicolumn{1}{c|}{\text{(j1)}} & \dt^B & 0 & -\dt^B \\* 
\hhline{----~~----}\nopagebreak
0 & \frac{23}{12} & -\frac{1}{2} & 0 &&& 
\dt^B & \frac{23}{12} & 0 & -\frac{5}{36} \\* 
-\dt^A & 0 & -1 & \frac{5}{12} &&& 
0 & 0 & -\frac{4}{3} & \frac{5}{12} \\* 
-2\dt^A & 0 & \frac{1}{6} & 0 &&& 
-\dt^A & 0 & 0 & \frac{5}{36} 
\vspace{\vgap} \\
\multicolumn{11}{c}{\textbf{Transitioning to GTS by increasing a step size}} \\*
\multicolumn{1}{c|}{\text{(k0)}} & 0 & -\dt^B & -2\dt^B &&& 
\multicolumn{1}{c|}{\text{(k1)}} & \dt^A & 0 & -\dt^B \\* 
\hhline{----~~----}\nopagebreak
0 & \frac{19}{6} & -\frac{5}{4} & 0 &&& 
\dt^A & \frac{37}{18} & 0 & -\frac{5}{36} \\* 
-\dt^A & 0 & -\frac{5}{2} & \frac{7}{6} &&& 
0 & 0 & -\frac{13}{6} & \frac{5}{6} \\* 
-2\dt^A & 0 & \frac{5}{12} & 0 &&& 
-\dt^A & 0 & 0 & \frac{5}{12} 
\end{longtable}

\subsection{Order 4}
\renewcommand\hgap{5pt}

\hspace{0pt} 

\begin{longtable}{s|sssss@{}p{\hgap}@{}s|ssss@{\hspace{-2.5pt}}s}
\multicolumn{13}{c}{\textbf{LTS $2:1$ rule}} \\*
\multicolumn{1}{c|}{\text{(a)}} & \dt^B & 0 & -\dt^B & -2\dt^B & -3\dt^B && 
\multicolumn{1}{c|}{\text{(b)}} & 0 & -\dt^B & -2\dt^B & -3\dt^B \\* 
\hhline{------~-----}\nopagebreak
0 & \frac{1925}{768} & -\frac{1}{12} & -\frac{55}{384} & 0 & \frac{3}{256} && 
0 & \frac{55}{24} & -\frac{295}{384} & 0 & \frac{3}{128} \\* 
-\dt^A & -\frac{1925}{768} & 0 & -\frac{55}{128} & \frac{7}{12} & -\frac{27}{256} && 
-\dt^A & 0 & -\frac{295}{128} & \frac{37}{24} & -\frac{27}{128} \\* 
-2\dt^A & \frac{385}{256} & 0 & \frac{55}{384} & 0 & -\frac{27}{256} && 
-2\dt^A & 0 & \frac{295}{384} & 0 & -\frac{27}{128} \\* 
-3\dt^A & -\frac{275}{768} & 0 & -\frac{11}{384} & 0 & \frac{3}{256} && 
-3\dt^A & 0 & -\frac{59}{384} & 0 & \frac{3}{128} 
\vspace{\vgap} \\
\multicolumn{1}{c|}{\text{(c)}} & \dt^B & 0 & -\dt^B & -2\dt^B \\* 
\hhline{-----}\nopagebreak
0 & \frac{1925}{384} & -\frac{59}{24} & \frac{185}{384} & 0 \\* 
-\dt^A & -\frac{1925}{384} & 0 & \frac{185}{128} & -\frac{3}{8} \\* 
-2\dt^A & \frac{385}{128} & 0 & -\frac{185}{384} & 0 \\* 
-3\dt^A & -\frac{275}{384} & 0 & \frac{37}{384} & 0 
\vspace{\vgap} \\
\multicolumn{13}{c}{\textbf{Transition to LTS by decreasing a step size}} \\*
\multicolumn{1}{c|}{\text{(d0)}} & \dt^B & 0 & -\dt^A & -2\dt^A & -3\dt^A && 
\multicolumn{1}{c|}{\text{(e0)}} & 0 & -\dt^A & -2\dt^A & -3\dt^A \\* 
\hhline{------~-----}\nopagebreak
0 & \frac{833}{384} & \frac{47}{384} & 0 & 0 & 0 && 
0 & \frac{99}{64} & 0 & 0 & 0 \\* 
-\dt^A & -\frac{833}{384} & 0 & -\frac{37}{128} & 0 & 0 && 
-\dt^A & 0 & -\frac{187}{192} & 0 & 0 \\* 
-2\dt^A & \frac{833}{640} & 0 & 0 & \frac{461}{1920} & 0 && 
-2\dt^A & 0 & 0 & \frac{107}{192} & 0 \\* 
-3\dt^A & -\frac{119}{384} & 0 & 0 & 0 & -\frac{25}{384} && 
-3\dt^A & 0 & 0 & 0 & -\frac{25}{192} 
\vspace{\vgap} \\
\multicolumn{1}{c|}{\text{(f0)}} & \dt^B & 0 & -\dt^A & -2\dt^A &&& 
\multicolumn{1}{c|}{\text{(d1)}} & 3\dt^B & 2\dt^B & \dt^B & 0 & -\dt^A \\* 
\hhline{-----~~------}\nopagebreak
0 & \frac{833}{192} & -\frac{125}{96} & 0 & 0 &&& 
\dt^A & \frac{1925}{768} & -\frac{25}{192} & -\frac{65}{768} & 0 & 0 \\* 
-\dt^A & -\frac{833}{192} & 0 & \frac{19}{48} & 0 &&& 
0 & -\frac{1925}{768} & 0 & -\frac{65}{256} & \frac{29}{96} & 0 \\* 
-2\dt^A & \frac{833}{320} & 0 & 0 & -\frac{37}{480} &&& 
-\dt^A & \frac{385}{256} & 0 & \frac{65}{768} & 0 & -\frac{3}{64} \\* 
-3\dt^A & -\frac{119}{192} & 0 & 0 & 0 &&& 
-2\dt^A & -\frac{275}{768} & 0 & -\frac{13}{768} & 0 & 0 
\vspace{\vgap} \\
\multicolumn{1}{c|}{\text{(e1)}} & 2\dt^B & \dt^B & 0 & -\dt^A &&& 
\multicolumn{1}{c|}{\text{(f1)}} & 3\dt^B & 2\dt^B & \dt^B & 0 \\* 
\hhline{-----~~-----}\nopagebreak
\dt^A & \frac{211}{96} & -\frac{125}{192} & 0 & 0 &&& 
\dt^A & \frac{1925}{384} & -\frac{59}{24} & \frac{185}{384} & 0 \\* 
0 & 0 & -\frac{125}{64} & \frac{47}{48} & 0 &&& 
0 & -\frac{1925}{384} & 0 & \frac{185}{128} & -\frac{3}{8} \\* 
-\dt^A & 0 & \frac{125}{192} & 0 & -\frac{3}{32} &&& 
-\dt^A & \frac{385}{128} & 0 & -\frac{185}{384} & 0 \\* 
-2\dt^A & 0 & -\frac{25}{192} & 0 & 0 &&& 
-2\dt^A & -\frac{275}{384} & 0 & \frac{37}{384} & 0 
\vspace{\vgap} \\
\multicolumn{13}{c}{\textbf{Transition to LTS by increasing a step size}} \\*
\multicolumn{1}{c|}{\text{(g0)}} & \dt^B & 0 & -\dt^B & -2\dt^B & -3\dt^B && 
\multicolumn{1}{c|}{\text{(h0)}} & 0 & -\dt^B & -2\dt^B & -3\dt^B \\* 
\hhline{------~-----}\nopagebreak
0 & \frac{55}{12} & -\frac{1}{12} & 0 & 0 & 0 && 
0 & \frac{55}{24} & 0 & 0 & 0 \\* 
-\dt^B & -\frac{55}{8} & 0 & -\frac{11}{24} & 0 & 0 && 
-\dt^B & 0 & -\frac{59}{24} & 0 & 0 \\* 
-2\dt^B & \frac{55}{12} & 0 & 0 & \frac{7}{12} & 0 && 
-2\dt^B & 0 & 0 & \frac{37}{24} & 0 \\* 
-3\dt^B & -\frac{55}{48} & 0 & 0 & 0 & -\frac{3}{16} && 
-3\dt^B & 0 & 0 & 0 & -\frac{3}{8} 
\vspace{\vgap} \\
\multicolumn{1}{c|}{\text{(i0)}} & \dt^B & 0 & -\dt^B & -2\dt^B &&& 
\multicolumn{1}{c|}{\text{(g1)}} & 3\dt^B & 2\dt^B & \dt^B & 0 & -\dt^B \\* 
\hhline{-----~~------}\nopagebreak
0 & \frac{55}{6} & -\frac{59}{24} & 0 & 0 &&& 
\dt^A & \frac{275}{96} & -\frac{1}{12} & -\frac{11}{96} & 0 & 0 \\* 
-\dt^B & -\frac{55}{4} & 0 & \frac{37}{24} & 0 &&& 
0 & -\frac{275}{48} & 0 & -\frac{11}{16} & \frac{7}{12} & 0 \\* 
-2\dt^B & \frac{55}{6} & 0 & 0 & -\frac{3}{8} &&& 
-\dt^B & \frac{275}{48} & 0 & \frac{11}{24} & 0 & -\frac{3}{16} \\* 
-3\dt^B & -\frac{55}{24} & 0 & 0 & 0 &&& 
-2\dt^B & -\frac{55}{32} & 0 & -\frac{11}{96} & 0 & 0 
\vspace{\vgap} \\
\multicolumn{1}{c|}{\text{(h1)}} & 2\dt^B & \dt^B & 0 & -\dt^B &&& 
\multicolumn{1}{c|}{\text{(i1)}} & 3\dt^B & 2\dt^B & \dt^B & 0 \\* 
\hhline{-----~~-----}\nopagebreak
\dt^A & \frac{55}{24} & -\frac{59}{96} & 0 & 0 &&& 
\dt^A & \frac{275}{48} & -\frac{59}{24} & \frac{37}{96} & 0 \\* 
0 & 0 & -\frac{59}{16} & \frac{37}{24} & 0 &&& 
0 & -\frac{275}{24} & 0 & \frac{37}{16} & -\frac{3}{8} \\* 
-\dt^B & 0 & \frac{59}{24} & 0 & -\frac{3}{8} &&& 
-\dt^B & \frac{275}{24} & 0 & -\frac{37}{24} & 0 \\* 
-2\dt^B & 0 & -\frac{59}{96} & 0 & 0 &&& 
-2\dt^B & -\frac{55}{16} & 0 & \frac{37}{96} & 0 
\vspace{\vgap} \\
\multicolumn{1}{c|}{\text{(g2)}} & 5\dt^B & 4\dt^B & 3\dt^B & 2\dt^B & \dt^B && 
\multicolumn{1}{c|}{\text{(h2)}} & 4\dt^B & 3\dt^B & 2\dt^B & \dt^B \\* 
\hhline{------~-----}\nopagebreak
2\dt^A & \frac{165}{64} & -\frac{1}{12} & -\frac{11}{80} & 0 & \frac{3}{320} && 
2\dt^A & \frac{55}{24} & -\frac{59}{80} & 0 & \frac{3}{160} \\* 
\dt^A & -\frac{275}{96} & 0 & -\frac{11}{24} & \frac{7}{12} & -\frac{3}{32} && 
\dt^A & 0 & -\frac{59}{24} & \frac{37}{24} & -\frac{3}{16} \\* 
0 & \frac{165}{64} & 0 & \frac{11}{48} & 0 & -\frac{9}{64} && 
0 & 0 & \frac{59}{48} & 0 & -\frac{9}{32} \\* 
-\dt^B & -\frac{55}{48} & 0 & -\frac{11}{120} & 0 & \frac{3}{80} && 
-\dt^B & 0 & -\frac{59}{120} & 0 & \frac{3}{40} 
\vspace{\vgap} \\
\multicolumn{1}{c|}{\text{(i2)}} & 5\dt^B & 4\dt^B & 3\dt^B & 2\dt^B \\* 
\hhline{-----}\nopagebreak
2\dt^A & \frac{165}{32} & -\frac{59}{24} & \frac{37}{80} & 0 \\* 
\dt^A & -\frac{275}{48} & 0 & \frac{37}{24} & -\frac{3}{8} \\* 
0 & \frac{165}{32} & 0 & -\frac{37}{48} & 0 \\* 
-\dt^B & -\frac{55}{24} & 0 & \frac{37}{120} & 0 
\vspace{\vgap} \\
\multicolumn{13}{c}{\textbf{Transitioning to GTS by decreasing a step size}} \\*
\multicolumn{1}{c|}{\text{(j0)}} & 0 & -\dt^B & -2\dt^B & -3\dt^B &&& 
\multicolumn{1}{c|}{\text{(j1)}} & \dt^B & 0 & -\dt^B & -2\dt^B \\* 
\hhline{-----~~-----}\nopagebreak
0 & \frac{55}{24} & -\frac{295}{384} & 0 & \frac{3}{128} &&& 
\dt^B & \frac{55}{24} & 0 & -\frac{37}{120} & 0 \\* 
-\dt^A & 0 & -\frac{295}{128} & \frac{37}{24} & -\frac{27}{128} &&& 
0 & 0 & -\frac{59}{24} & \frac{37}{32} & 0 \\* 
-2\dt^A & 0 & \frac{295}{384} & 0 & -\frac{27}{128} &&& 
-\dt^A & 0 & 0 & \frac{37}{48} & -\frac{3}{8} \\* 
-3\dt^A & 0 & -\frac{59}{384} & 0 & \frac{3}{128} &&& 
-2\dt^A & 0 & 0 & -\frac{37}{480} & 0 
\vspace{\vgap} \\
\multicolumn{1}{c|}{\text{(j2)}} & 2\dt^B & \dt^B & 0 & -\dt^B \\* 
\hhline{-----}\nopagebreak
2\dt^B & \frac{55}{24} & 0 & 0 & -\frac{3}{32} \\* 
\dt^B & 0 & -\frac{59}{24} & 0 & \frac{3}{8} \\* 
0 & 0 & 0 & \frac{37}{24} & -\frac{9}{16} \\* 
-\dt^A & 0 & 0 & 0 & -\frac{3}{32} 
\vspace{\vgap} \\
\multicolumn{13}{c}{\textbf{Transitioning to GTS by increasing a step size}} \\*
\multicolumn{1}{c|}{\text{(k0)}} & 0 & -\dt^B & -2\dt^B & -3\dt^B &&& 
\multicolumn{1}{c|}{\text{(k1)}} & \dt^A & 0 & -\dt^B & -2\dt^B \\* 
\hhline{-----~~-----}\nopagebreak
0 & \frac{9}{2} & -\frac{55}{24} & 0 & \frac{1}{12} &&& 
\dt^A & \frac{8}{3} & 0 & -\frac{3}{8} & 0 \\* 
-\dt^A & 0 & -\frac{55}{8} & \frac{31}{6} & -\frac{3}{4} &&& 
0 & 0 & -\frac{35}{6} & \frac{27}{8} & 0 \\* 
-2\dt^A & 0 & \frac{55}{24} & 0 & -\frac{3}{4} &&& 
-\dt^A & 0 & 0 & \frac{27}{8} & -\frac{11}{6} \\* 
-3\dt^A & 0 & -\frac{11}{24} & 0 & \frac{1}{12} &&& 
-2\dt^A & 0 & 0 & -\frac{3}{8} & 0 
\vspace{\vgap} \\
\multicolumn{1}{c|}{\text{(k2)}} & 2\dt^A & \dt^A & 0 & -\dt^B \\* 
\hhline{-----}\nopagebreak
2\dt^A & \frac{71}{30} & 0 & 0 & -\frac{3}{40} \\* 
\dt^A & 0 & -\frac{17}{6} & 0 & \frac{3}{8} \\* 
0 & 0 & 0 & \frac{8}{3} & -\frac{9}{8} \\* 
-\dt^A & 0 & 0 & 0 & -\frac{3}{8} 
\end{longtable}

\endgroup